\documentclass[]{elsarticle}

\usepackage[margin=2.5cm]{geometry}
\usepackage{framed,multirow}
\usepackage{mathptmx}

\usepackage{amssymb}
\usepackage{latexsym}

\usepackage{url}
\usepackage{xcolor}
\definecolor{newcolor}{rgb}{.8,.349,.1}

\usepackage{tcolorbox}
\tcbuselibrary{theorems}
\usepackage{xcolor} 
\usepackage[noend]{algpseudocode}
\usepackage{array}

\usepackage{amsfonts}
\usepackage{graphicx}
\usepackage{bm,mathrsfs}
\usepackage{epstopdf}
\usepackage{mathtools}
\usepackage{tikz,pgfplots}
\usetikzlibrary{snakes,arrows,shapes,trees}
\usetikzlibrary{arrows.meta,calc,positioning}
\usepackage{amsmath,amsthm}
\usepackage{amsopn}
\usepackage{listings}
\usepackage{adjustbox}
\usepackage{longtable}
\usepackage{multirow}
\usepackage{hyperref}
\usepackage{pgfplots}
\usepackage{pgfplotstable}
\usepackage{grffile}
\usetikzlibrary{arrows,shapes,plotmarks}
\usepgfplotslibrary{groupplots}
\usetikzlibrary{matrix}
\usepackage{siunitx}
\usepackage[switch]{lineno}
\usepackage{enumitem}

\usepackage{titlesec}
\newcommand{\cref}[2]{\hyperref[#2]{#1~\ref*{#2}}}
\newcommand{\figref}[1]{\hyperref[#1]{Fig.~\ref*{#1}}}
\newcommand{\secref}[1]{\hyperref[#1]{Sec.~\ref*{#1}}}
\newcommand{\tabref}[1]{\hyperref[#1]{Tab.~\ref*{#1}}}
\newcommand{\eqnref}[1]{\hyperref[#1]{Eq.~(\ref*{#1})}}
\newcommand{\Algref}[1]{\hyperref[#1]{Algorithm~\ref*{#1}}}
\newcommand{\pd}[2]{\frac{\partial #1}{\partial #2}} 
 
\hypersetup{
	colorlinks,
	linkcolor={blue!50!black},
	citecolor={blue!50!black},
	urlcolor={blue!80!black},
	anchorcolor = {blue!80!black},
	filecolor = {blue!80!black},
	menucolor = {blue!80!black},
	runcolor = {blue!80!black}
}

\pgfplotsset{compat=1.8}
\usepackage{soul}
\usepackage{rotating}
\usepackage{url}
\usepackage{algorithm}

\usepackage{enumitem}
\usepackage{subcaption}

\newcommand{\Vector}[1]{\mathbf{#1}}

\newcommand{\Tensor}[1]{\underline{\underline{\mathbf{#1}}}}

\algnewcommand{\LeftComment}[1]{\Statex \(\triangleright\) #1}

\pgfplotsset{
compat=1.8,
legend image code/.code={
\draw[mark repeat=2,mark phase=2]
plot coordinates {
(0cm,0cm)
(0.15cm,0cm)        
(0.3cm,0cm)         
};%
}
}

\newcommand{\norm}[1]{\left\lVert#1\right\rVert}
\definecolor{cpu3}{HTML}{F44336}
\definecolor{cpu4}{HTML}{2196F3}
\definecolor{cpu1}{HTML}{4CAF50}
\definecolor{cpu2}{HTML}{FFC107}
\definecolor{gpu3}{HTML}{EF9A9A}
\definecolor{gpu4}{HTML}{90CAF9}
\definecolor{gpu1}{HTML}{A5D6A7}
\definecolor{gpu2}{HTML}{FFE082}

\definecolor{cpu5}{HTML}{9932CC}

\definecolor{sq_b1}{RGB}{37,52,148}
\definecolor{sq_b2}{RGB}{44,127,184}
\definecolor{sq_b3}{RGB}{65,182,196}
\definecolor{sq_b4}{RGB}{127,205,187}
\definecolor{sq_b5}{RGB}{199,233,180}
\definecolor{sq_b6}{RGB}{255,255,204}

\definecolor{sq_r1}{RGB}{189,0,38}
\definecolor{sq_r2}{RGB}{240,59,32}
\definecolor{sq_r3}{RGB}{253,141,60}
\definecolor{sq_r4}{RGB}{254,178,76}
\definecolor{sq_r5}{RGB}{254,217,118}
\definecolor{sq_r6}{RGB}{255,255,178}

\definecolor{sq_g1}{RGB}{0,104,55}
\definecolor{sq_g2}{RGB}{49,163,84}
\definecolor{sq_g3}{RGB}{120,198,121}
\definecolor{sq_g4}{RGB}{173,221,142}
\definecolor{sq_g5}{RGB}{217,240,163}
\definecolor{sq_g6}{RGB}{255,255,204}

\definecolor{div_c1}{RGB}{230,171,2}
\definecolor{div_c2}{RGB}{102,166,30}
\definecolor{div_c3}{RGB}{231,41,138}
\definecolor{div_c4}{RGB}{117,112,179}
\definecolor{div_c5}{RGB}{217,95,2}
\definecolor{div_c6}{RGB}{27,158,119}
\definecolor{div_c7}{RGB}{215,48,39}

\definecolor{div_d1}{RGB}{215,25,28}
\definecolor{div_d2}{RGB}{253,174,97}
\definecolor{div_d3}{RGB}{255,255,191}
\definecolor{div_d4}{RGB}{171,217,233}
\definecolor{div_d5}{RGB}{44,123,182}
\definecolor{ao}{RGB}{0.0, 128, 0.0}
\usepackage{dirtytalk}

\newcommand{\added}[1]{\textcolor{black}{#1}}
\usepackage[normalem]{ulem}

\usepgfplotslibrary{fillbetween}
 
 \usepackage{lineno}
 \renewcommand{\vec}[1]{\mathbf{#1}}

 \usetikzlibrary{external}
\begin{document}


\begin{frontmatter}

\title{Field conserving adaptive mesh refinement (AMR) scheme on massively parallel adaptive octree meshes}%

\author[ISU]{Kumar Saurabh}
\ead{maksbh@iastate.edu}
\author[Stanford]{Makrand A. Khanwale\texorpdfstring{\corref{cor}}{}}
\ead{khanwale@stanford.edu}
\author[Utah,Riken]{Masado Ishii}
\ead{masado.ishii@riken.jp}
\author[Utah,Tufts]{Hari Sundar}
\ead{hari.sundar@tufts.edu}
\author[ISU]{Baskar Ganapathysubramanian\texorpdfstring{\corref{cor}}{}}
\ead{baskarg@iastate.edu}
\cortext[cor]{Corresponding authors}
\address[ISU]{Department of Mechanical Engineering, Iowa State University, Ames, IA}
\address[Stanford]{Department of Mechanical Engineering, Stanford University, Stanford, CA}
\address[Utah]{Kahlert School of Computing, University of Utah, Salt Lake City, Utah}
\address[Riken]{RIKEN Center for Computational Science, Hyogo, Japan}
\address[Tufts]{Department of Computer Science, Tufts University, Medford, MA}
\begin{abstract}


Adaptive mesh refinement (AMR) is widely used to efficiently resolve localized features in time-dependent partial differential equations (PDEs) by selectively refining and coarsening the mesh. However, in long-horizon simulations, repeated intergrid interpolations can introduce systematic drift in conserved quantities, especially for variational discretizations with continuous basis functions. While interpolation from parent-to-child during refinement in continuous Galerkin (CG) discretizations is naturally conservative, the standard injection-based child-to-parent coarsening interpolation is generally not.

We propose a simple, scalable \textit{field-conserving coarsening} operator for parallel, \emph{octree}-based AMR. The method enforces discrete global conservation during coarsening by first computing field conserving coarse-element values at quadrature points and then recovering coarse nodal degrees of freedom via an $L^2$ projection (mass-matrix solve), which simultaneously controls the $L_2$ error. 
We evaluate the approach on mass-conserving phase-field models, including the Cahn--Hilliard and Cahn--Hilliard--Navier--Stokes systems, and compare against injection in terms of conservation error, solution quality, and computational cost. 
\end{abstract}


\end{frontmatter}


\section{Introduction}

Mass conservation is a fundamental requirement in many multiphysics simulations, and it becomes especially critical in long-horizon computations where small per-step defects can accumulate into noticeable drift. Numerous studies have aimed to derive governing equations that inherently conserve mass, such as those used for simulating phase separation phenomena like the Cahn--Hilliard equation \citep{elliott1989cahn,gameiro2005evolution,speck2014effective,blowey1991cahn}, or for interface tracking methods including Cahn--Hilliard--Navier--Stokes \citep{han2015second,diegel2017convergence,khanwale2020simulating,boyer2010cahn} and conservative level-set approaches \citep{olsson2005conservative,olsson2007conservative,van2005mass}. Additionally, several works have focused on developing accurate numerical schemes that successfully simulate these PDEs and have demonstrated that the underlying invariants (e.g., mass) are conserved to numerical precision on a fixed mesh \citep{khanwale2020simulating,Khanwale2023projection,khanwale2022fully,wodo2011computationally,tornberg2000finite,marchandise2006stabilized}.

In practical simulations of these PDEs, adaptive mesh refinement (AMR) is often essential: it concentrates resolution near evolving interfaces or sharp gradients while maintaining tractable cost. Among AMR strategies, \emph{balanced adaptive octree meshes} are widely used due to their algorithmic simplicity, strict locality of refinement/coarsening, and scalability on large core counts \citep{tu2005scalable,ishii2019solving,fernando2018massively,sundar2008bottom,bastian2008generic,greaves1999hierarchical,bader2012space,BursteddeWilcoxGhattas11,popinet2003gerris,saurabh2022scalable,saurabh2021industrial,saurabh2021scalable}. In octree AMR, adaptation proceeds through two fundamental operations: \emph{refinement}, which subdivides selected octants, and \emph{coarsening}, which merges sibling octants. Both operations require transferring discrete fields between meshes of different resolutions. The transfer of a field from one grid ($\mathcal{M}_1$) to another ($\mathcal{M}_2$) is commonly referred to as \textit{intergrid transfer} and can be written as
\begin{equation}
    \Tensor{T}(\phi_{\mathcal{M}_1}) = \phi_{\mathcal{M}_2},
\end{equation}
where $\phi_{\mathcal{M}_1}$ and $\phi_{\mathcal{M}_2}$ denote the discrete field representations on $\mathcal{M}_1$ and $\mathcal{M}_2$, respectively. For a conserved scalar $\phi$, a conservative transfer satisfies
\begin{equation}
    \bigg[\int_{\Omega} \phi \; dV\bigg]_{\mathcal{M}_1} = \bigg[\int_{\Omega} \phi \; dV\bigg]_{\mathcal{M}_2}.
    \label{Eq:intergrid_transfer}
\end{equation}
A stronger requirement (which our approach satisfies) is for Eq.~\ref{Eq:intergrid_transfer} to be valid for any arbitrary sub-domain $\Omega_k \subset \Omega$.

Although refinement and coarsening are both common in AMR, they play fundamentally different roles with respect to conservation in continuous Galerkin (CG) discretizations. In CG, refinement transfers can be performed using basis-function interpolation from parent to child and are naturally conservative for integral quantities of the represented field. In contrast, the standard injection-based coarsening transfer (child to parent) removes fine-grid degrees of freedom that do not coincide with coarse-grid nodes, which can introduce systematic drift in global invariants. Over long time horizons, precisely the regime of interest for many multiphysics problems, this drift can dominate the error budget even when the underlying PDE and time integrator are conservative.

This work targets this failure mode in the specific and widely used setting of \emph{parallel octree-based AMR with single-level adaptation}. We show that, for CG finite element discretizations on octree meshes, non-conservative coarsening is the dominant source of global mass drift, and we propose a simple, scalable \emph{field-conserving coarsening operator} that guarantees discrete global conservation across AMR cycles. While our demonstrations and numerical examples use linear and quadratic Lagrange elements for clarity, the proposed coarsening procedure applies to \emph{any} CG formulation with Lagrange basis functions of arbitrary polynomial order (e.g., $Q_p$/$P_p$, $p\ge 1$), provided that the transfer and projection use quadrature consistent with the chosen order.

Octree AMR proceeds in two mutually exclusive stages:
\begin{enumerate}
\item \textbf{Refine:} each element in $\mathcal{M}_2$ is either at a finer level or at the same level as in $\mathcal{M}_1$.
\item \textbf{Coarsen:} each element in $\mathcal{M}_2$ is either at a coarser level or at the same level as in $\mathcal{M}_1$.
\end{enumerate}
During the \textbf{refine} stage, each octant is flagged as either \texttt{REFINE} or \texttt{NO\_CHANGE}, while during the \textbf{coarsen} stage it is flagged as \texttt{COARSEN} or \texttt{NO\_CHANGE}. We restrict refinement and coarsening to a single level per adaptation cycle (splitting/merging $2/4/8$ octants in 1D/2D/3D), which is standard in octree AMR implementations \citep{sundar2008bottom,fernando2018massively,ishii2019solving,saurabh2022scalable}.

\subsection{Conservation properties of octree AMR intergrid transfer}
\label{subsec:conservation_transfer}

We briefly analyze the refine and coarsen stages to motivate why conservative \emph{coarsening} is essential for long-time simulations.

\textbf{Refinement is conservative (parent-to-child transfer):} In the refine stage, the intergrid operator $\Tensor{T}_R$ interpolates values from a parent element to its children (also referred to as parent-to-child, \texttt{P2C}) using the basis functions of the CG discretization.\footnote{In octree AMR, the coarser element is called the \emph{parent} of the finer element and each refined element is called a \emph{child}.} For linear elements, this implies that the midpoint value introduced by refinement is consistent with the linear representation, and hence the integral of the represented field is preserved.

\begin{figure}[h]
	\centering
\begin{subfigure}[b]{1.0\linewidth}
	\centering
\begin{tikzpicture}
\begin{axis}[
          width=1.0\linewidth, 
          height = 0.5\linewidth,
            xlabel={$x$},
    ylabel={$g(x)$},
    axis x line=bottom,
    axis y line=left,
	legend style={at={(0.50,0.35)},anchor= north,nodes={scale=0.95, transform shape}},
]
     \addplot[only marks,color=blue]
        table[x expr={\thisrow{x}))},y expr=(\thisrow{val},col sep=space]{Data/Curve1D/fineMesh.txt};
     \addplot +[mark=square,color = red,  thick] 
        table[x expr={\thisrow{x}))},y expr=(\thisrow{val},col sep=space]{Data/Curve1D/refinement.txt};

        \legend{\small{Original mesh},\small{Refined mesh}}
  \end{axis}
  \end{tikzpicture}
  \end{subfigure}
    \caption{\textbf{Refinement:} Figure representing the refinement case. The red marker represents the extra points that are added on the refined mesh. The value of 
    	$\bigg[\int_{\Omega} g(x) \; d\Omega \bigg]_{\mathcal{M}_O}$ = 
    	$\bigg[\int_{\Omega} g(x) \; d\Omega \bigg]_{\mathcal{M}_R}$ = 10.6284 , where $\mathcal{M}_O$ is the original mesh and $\mathcal{M}_R$ is the refined mesh.
}
    \label{fig:Refinement1D}
\end{figure}
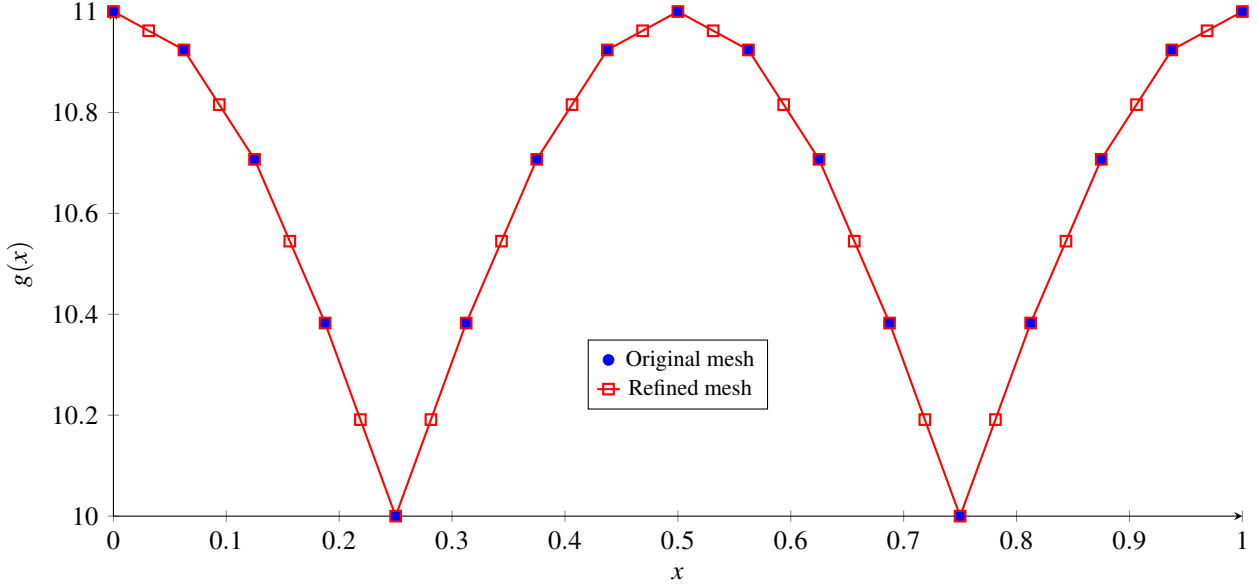

To illustrate, consider a 1D linear function $g(x)=\alpha x + \beta$. The area under the curve over $[x_1,x_2]$ is
\begin{equation}
    A = \int_{x_1}^{x_2} g(x)\,dx
      = \frac{\alpha (x_2^2 - x_1^2)}{2} + \beta (x_2 - x_1).
    \label{eq:refine_1}
\end{equation}
Refinement in 1D splits the parent element by introducing $x'=(x_1+x_2)/2$. Since we interpolate using the linear basis, $g(x')=\alpha(x_1+x_2)/2+\beta$, and
\begin{equation*}
\begin{split}
    A &= \int_{x_1}^{x_2} g(x)\,dx
      = \int_{x_1}^{x'} g(x)\,dx + \int_{x'}^{x_2} g(x)\,dx \\
      &= \frac{\alpha (x_2^2 - x_1^2)}{2} + \beta (x_2 - x_1),
\end{split}
\end{equation*}
which is identical to \eqnref{eq:refine_1}. Thus, the refinement stage preserves the integral of the field variable. \figref{fig:Refinement1D} illustrates this behavior for $g(x)=|\cos(2\pi x)|+10$. Algorithmically, refine-stage transfer is local: new nodal values are computed from parent degrees of freedom. We refer to \citet{saurabh2022scalable} for implementation details.

\textbf{Injection-based coarsening is not conservative (child-to-parent transfer):} The coarsening transfer operator $\Tensor{T}_C$ is more subtle. The common state-of-the-art procedure uses an injection-based approach, where child degrees of freedom that do not coincide with parent nodes are eliminated from the representation (also referred to as child-to-parent, \texttt{C2P}, transfer). This procedure generally does \emph{not} conserve integral quantities of the represented field.

\begin{figure}[h]
\begin{subfigure}[b]{1.0\linewidth}
	\centering
\begin{tikzpicture}
\begin{axis}[
          width=1.0\linewidth, 
          height = 0.5\linewidth,
            xlabel={$x$},
    ylabel={$g(x)$},
    axis x line=bottom,
    axis y line=left,
	legend style={at={(0.50,0.35)},anchor= north,nodes={scale=0.95, transform shape}},
]
     \addplot[mark=o,color = blue,  thick,name path=A] 
        table[x expr={\thisrow{x}))},y expr=(\thisrow{val},col sep=space]{Data/Curve1D/fineMesh.txt};
     \addplot[mark=o,color = red,  thick,name path=B] 
        table[x expr={\thisrow{x}))},y expr=(\thisrow{val},col sep=space]{Data/Curve1D/injection.txt};
           \addplot[black!40] fill between[of=A and B];
           \legend{\small{Original mesh}, \small{Coarsened mesh}}
  \end{axis}
  \end{tikzpicture}
  \end{subfigure}
  \caption{\textbf{Coarsening by Injection}: Figure demonstrating the coarsening procedure by injection. The gray shaded region shows the area corresponding to the difference between two mesh. 	$\bigg[\int_{\Omega} g(x) \; d\Omega \bigg]_{\mathcal{M}_O}$ = 10.6284, whereas 
  	$\bigg[\int_{\Omega} g(x) \; d\Omega \bigg]_{\mathcal{M}_C}$ = 10.6036 , where $\mathcal{M}_O$ is the original mesh and $\mathcal{M}_C$ is the coarsened mesh.}
  \label{fig:Coarsening1D}
\end{figure}
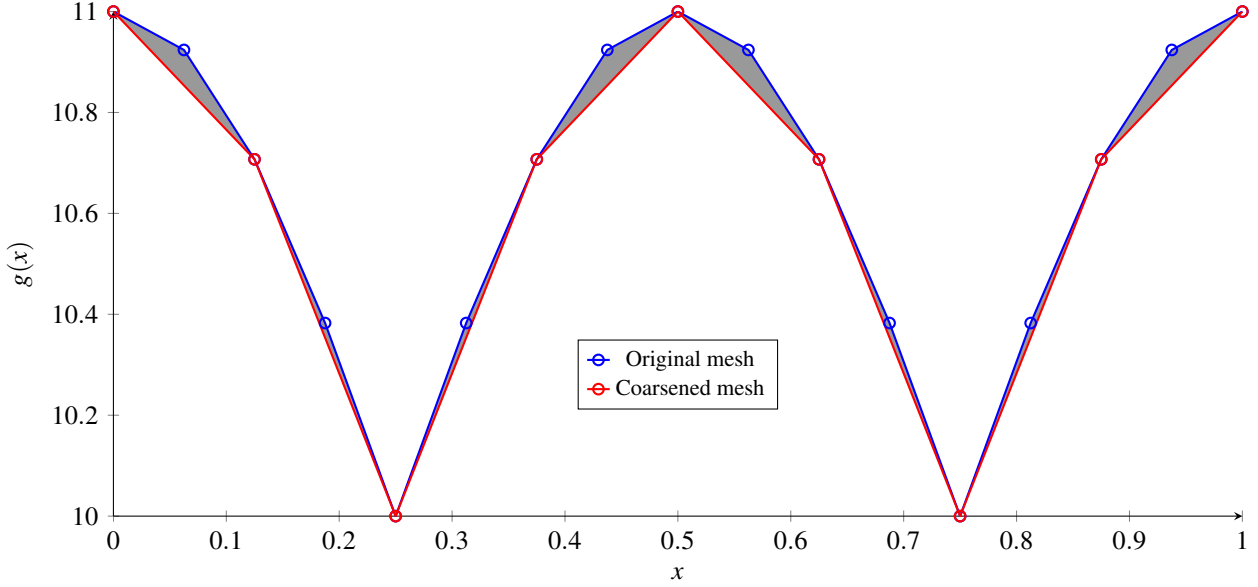

To see this, consider a piecewise linear function over $[x_1,x_3]$ with a breakpoint at $x_2$:
\begin{equation}
 g(x) =
 \begin{cases}
  \alpha_1 x + \beta_1  & x \in[x_1,x_2], \\
  \alpha_2 x + \beta_2  & x \in[x_2,x_3].
 \end{cases}
\end{equation}
Coarsening merges the two child elements into one parent element by removing the interior node $x_2$. Injection connects $(x_1,g(x_1))$ to $(x_3,g(x_3))$ with a new line
\begin{equation*}
    g_C(x) = \alpha_C x + \beta_C, \qquad x\in[x_1,x_3],
\end{equation*}
such that $g_C(x_1)=g(x_1)$ and $g_C(x_3)=g(x_3)$, where
\begin{equation*}
\begin{split}
    \alpha_C & = \frac{\alpha_2 x_3 + \beta_2 - \alpha_1 x_1 - \beta_1}{x_3 - x_1},\\
    \beta_C  & = \frac{(\alpha_1 - \alpha_2)x_1x_3 + \beta_1 x_3 - \beta_2 x_1}{x_3 - x_1}.
\end{split}
\end{equation*}
In general, $\int_{x_1}^{x_3} g(x)\,dx \neq \int_{x_1}^{x_3} g_C(x)\,dx$, and equality holds only when $\alpha_1=\alpha_2$ and $\beta_1=\beta_2$ (i.e., when the two children lie on a single straight line). \figref{fig:Coarsening1D} demonstrates this effect for the same test function $g(x)=|\cos(2\pi x)|+10$, where the shaded region indicates the loss (or gain) in the integrated quantity caused by injection-based coarsening.

\textbf{Why CG requires a global approach:} A key distinction between continuous Galerkin (CG) and local discretization methods (e.g., finite volume or discontinuous Galerkin (DG)) is that CG basis functions have global support across element boundaries. In finite-volume or DG methods, mass can be redistributed locally during coarsening by conservatively averaging cell values or redistributing mass to neighboring elements within a local stencil. In contrast, CG's continuous inter-element representation means that altering nodal values to enforce conservation at coarsening sites necessarily affects the represented field across multiple elements. Consequently, restoring global conservation in CG cannot be achieved solely locally when degrees of freedom are eliminated during coarsening; instead, it requires a mechanism that involves all degrees of freedom (e.g., via a global constraint or projection).

\section{Numerical procedure}
\label{sec:Numerics}
We consider a scalar field $g(\vec{x})$ represented in a continuous Galerkin (CG) finite element space on a \emph{balanced adaptive octree mesh}. Our goal is to transfer this field from a \emph{fine} mesh to the mesh obtained after \emph{single-level coarsening} (merging $2/4/8$ sibling octants in 1D/2D/3D) while \emph{guaranteeing discrete global conservation} of the integral quantity (e.g., mass). The key difficulty is that, under standard injection-based child-to-parent transfer, degrees of freedom (DOFs) that do not coincide with the coarse mesh are discarded, which generally produces conservation drift over repeated AMR cycles. We therefore construct a conservative coarsening operator tailored to octree AMR that (i) conservatively transfers information to coarse-element quadrature points and (ii) reconstructs coarse nodal DOFs via an $L^2$ projection (mass-matrix solve).

Although we illustrate the method using linear Lagrange elements and common Gauss--Legendre stencils, the procedure applies to \emph{any} CG discretization with Lagrange basis functions of arbitrary polynomial degree ($P_p$/$Q_p$, $p\ge 1$), provided the quadrature rule used for transfer and projection is consistent with the chosen polynomial order.

\subsection{Discrete conservation constraint at quadrature points}
\label{subsec:disc_consv}

Let a single coarse (parent) element be formed by coarsening $N_{\mathrm{elem},f}$ fine (child) elements, where $N_{\mathrm{elem},f}=2/4/8$ in 1D/2D/3D, respectively. Let $\{(\vec{x}_{i,j}^f,w_{i}^f)\}_{i=1}^{N_{\mathrm{ip},f}}$ denote the quadrature points and weights on the $j$-th fine element, and let $\{(\vec{x}_{k}^c,w_{k}^c)\}_{k=1}^{N_{\mathrm{ip},c}}$ denote the quadrature points and weights on the resulting coarse element. Denote by $g^f(\vec{x}_{i,j}^f)$ the fine-mesh field evaluated at fine quadrature points, and by $g_c(\vec{x}_k^c)$ the (unknown) coarse-element quadrature values to be determined.

In the discrete setting, global conservation of the integral quantity over the parent element requires that the sum of the child contributions equals the coarse contribution (to numerical precision):
\begin{equation}
	\sum_{j=1}^{N_{\mathrm{elem},f}}\;\sum_{i=1}^{N_{\mathrm{ip},f}}
        w_i^f\, g^f\!\left(\vec{x}_{i,j}^f\right)
    \;=\;
    \sum_{k=1}^{N_{\mathrm{ip},c}}
        w_k^c\, g_c\!\left(\vec{x}_k^c\right).
	\label{eqn:mass_consv_disc}
\end{equation}
This essentially is a discrete statement of Eq.~\ref{Eq:intergrid_transfer} applied to a single coarse (parent) element.

\begin{figure}[H]
  \begin{subfigure}{0.45\textwidth}
    \begin{tikzpicture}[scale = 0.8]
    \draw (0,0) node[circle, fill, inner sep=2pt] {} -- (8,0) node[circle, fill, inner sep=2pt] {};
    \draw (4,0) node[circle, fill, inner sep=2pt] {};
    
    \draw (0.845,0) node[magenta, scale=2.0] {$\bm{\times}$};
    \draw (3.155,0) node[magenta, scale=2.0] {$\bm{\times}$};

     \draw (4.845,0) node[magenta, scale=2.0] {$\bm{\times}$};
    \draw (7.155,0) node[magenta, scale=2.0] {$\bm{\times}$}; 
    
    \draw (0,-1) node[circle, fill, inner sep=2pt] {} -- (8,-1) node[circle, fill, inner sep=2pt] {};

    \draw (1.691,-1) node[magenta, scale=2.0] {$\bm{\times}$};
    \draw (6.309,-1) node[magenta, scale=2.0] {$\bm{\times}$};
       \draw[<->, thick] (0,0.5) -- (4,0.5) node[midway, above] {$e_0^f$};
    \draw[<->, thick] (4,0.5) -- (8,0.5) node[midway, above] {$e_1^f$};
     \draw[<->, thick] (0,-1.5) -- (8,-1.5) node[midway, below] {$e_0^c$};
\end{tikzpicture}
\caption{Q1 element and quadrature points}
\label{fig:Q1_schematic}
  \end{subfigure}
\begin{subfigure}{0.45\textwidth}
    \begin{tikzpicture}[scale = 0.8]
    \draw (0,0) node[circle, fill, inner sep=2pt] {} -- (8,0) node[circle, fill, inner sep=2pt] {};

     \draw (0.451,0) node[magenta, scale=2.0] {$\bm{\times}$};
      \draw (2.0,0) node[magenta, scale=2.0] {$\bm{\times}$};
    \draw (3.549,0) node[magenta, scale=2.0] {$\bm{\times}$};
    \draw (4.451,0) node[magenta, scale=2.0] {$\bm{\times}$};
      \draw (6.0,0) node[magenta, scale=2.0] {$\bm{\times}$};
    \draw (7.549,0) node[magenta, scale=2.0] {$\bm{\times}$};
    \draw[<->, thick] (0,0.5) -- (4,0.5) node[midway, above] {$e_0^f$};
    \draw[<->, thick] (4,0.5) -- (8,0.5) node[midway, above] {$e_1^f$};
 \draw (4,0) node[circle, fill, inner sep=2pt] {};
    \draw (2,0) node[circle, fill, inner sep=2pt] {};
    \draw (6,0) node[circle, fill, inner sep=2pt] {};
    \draw (0,-1) node[circle, fill, inner sep=2pt] {} -- (8,-1) node[circle, fill, inner sep=2pt] {};
    \draw (4,-1) node[circle, fill, inner sep=2pt] {};
     \draw[<->, thick] (0,-1.5) -- (8,-1.5) node[midway, below] {$e_0^c$};
    
    \draw (0.902,-1) node[magenta, scale=2.0] {$\bm{\times}$};
    \draw (4,-1) node[magenta, scale=2.0] {$\bm{\times}$};
    \draw (7.098,-1) node[magenta, scale=2.0] {$\bm{\times}$};
    \node[circle, fill, inner sep=2pt] at (4,-1) {};
    \node[right] at (8.5,0.5) {Level};
    \node[right] at (8.5,0) {$l+1$};
    \node[right] at (8.5,-1) {$l$};
\end{tikzpicture}
\caption{Q2 element and quadrature points}
\label{fig:Q2_schematic}
  \end{subfigure}
  \centering
  \begin{tikzpicture}
    \node[magenta, scale=2.0] at (8.0,0) {$\bm{\times}$};
    \node[right] at (8.1,0) {Quadrature points};
     \node[circle, fill, inner sep=2pt] at (8.0,0.5) {};
    \node[right] at (8.1,0.5) {Nodal points};
  \end{tikzpicture}
\caption{Figure showing Gauss quadrature points and nodal points for (a) Q1 and (b) Q2 elements at two levels of refinement in 1D.}
\label{fig:schematic_1D}
\end{figure}
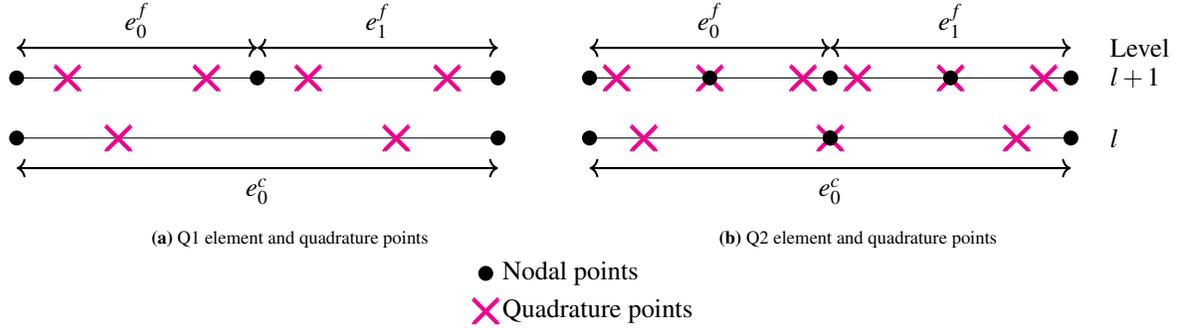

\textbf{Quadrature choice and mapping:} During the coarsening step, the left side of  the \eqref{eqn:mass_consv_disc} is fully determined by the fine-mesh representation. To compute the right side, we must choose a quadrature rule on the coarse element and determine a mapping between fine elements and coarse quadrature points. \figref{fig:schematic_1D} illustrates the location of standard Gauss-Legendre quadrature points and nodal points for $Q_1$ and $Q_2$ elements in 1D. We want to interpolate the values from the fine element (denoted by $e_0^f$ and $e_1^f$) quadrature points to the coarse element quadrature points. To achieve this goal, we construct a restriction operator by performing a local $L^2$ projection from the fine elements to the coarse element.

\tikzexternaldisable
\tikzexternalenable 

\subsection{1D local $L^2$ projection example}
Consider a 1D domain with two fine elements ($e_0^f$ and $e_1^f$) that are to be coarsened into a single coarse element ($e_0^c$). Let the parent reference element be $[-1,1]$ and the two child reference elements be $[-1,0]$ and $[0,1]$. We assume the polynomial space $Q_p$ with $p=1,2$ for both fine and coarse elements. We use the quadrature for numerical integration with $N_{\mathrm{ip},f}$ and $N_{\mathrm{ip},c}$ quadrature points for fine and coarse elements, respectively. \figref{fig:schematic_1D} shows the schematic for Gauss-Legendre quadrature points for $Q_1$ and $Q_2$ elements with $N_{\mathrm{ip},c}$ = $N_{\mathrm{ip},f}$ = $(p + 1)$.

\textbf{Coarse representation:}  We represent the coarse field on the parent element as:
\begin{equation}
g_c(x) = \sum_{j=1}^{N_{\mathrm{ip},c}} U_j N_j(x)
\label{eq:coarse_representation}
\end{equation}
where $N_j(x)$ are the Lagrange basis functions and $U_j$ are the nodal values at the coarse element nodes. We first consider the typical case of standard Gauss--Legendre quadrature with $N_{\mathrm{ip},c}=N_{\mathrm{ip},f}=p+1$ points.  Additionally, the polynomials are constructed such that the nodes are defined at the quadrature points: $N_j(\vec{x}_{k}^c) = \delta_{jk}$. Therefore, the coarse degrees of freedom are the values of the field at the quadrature points, i.e., $g_c(\vec{x}_{k,i}^c) = U_j$. We make this choice deliberately to obtain a diagonal mass matrix in the local $L^2$ projection defined below.

\textbf{Local $L^2$ projection definition:} Given the fine element values $g^f$ on the children, we want to compute the coarse coefficients $U_j$ such that the $L^2$ error between the fine element representation and the coarse element representation is minimized. This can be mathematically stated as:
\begin{equation}
  \label{eq:l2_min}
  \{U_j\}_{j=1}^{N_{\mathrm{ip},c}}
  = \arg\min_{\{U_j\}}
  \int_{K}\Big(g_f(x) - g_c(x)\Big)^2\,dx,
\end{equation}
where $K$ is the parent element, $g_f$ denotes the (piecewise) fine representation on the two children, and $g_c$ is the coarse representation on $K$.
The minimization in~\eqnref{eq:l2_min} can be posed as a Galerkin orthogonality condition as, 
\begin{equation}
  \label{eq:galerkin}
  \int_{K}\big(g_f(x) - g_c(x)\big)\,v(x)\,dx = 0
  \qquad \forall v\in Q_p(K).
\end{equation}
Choosing $v=N_i$ and using the discrete representation of $g_c$ from~\eqnref{eq:coarse_representation} yields the following local mass-matrix system
\begin{equation}
  \label{eq:mass_system}
  \sum_{j=1}^{N_{\mathrm{ip},c}} M_{ij}\,U_j = b_i,
  \qquad
  M_{ij}=\int_{K}N_i(x)N_j(x)\,dx,
  \qquad
  b_i=\int_{K} g_f(x)\,N_i(x)\,dx.
\end{equation}

To evaluate $b_i$ from child data, we split the parent into left and right children, $K=K_L\cup K_R$, and map each child reference coordinate $r\in[-1,1]$ to the parent coordinate $x\in[-1,1]$ via
\begin{equation}
  \label{eq:child_maps}
  x_L(r)=\tfrac12(r-1),\qquad x_R(r)=\tfrac12(r+1),\qquad dx=\tfrac12\,dr.
\end{equation}
 The right-hand side $b_i$ can be computed by summing the value over each child with nodes/weights $\{(r_q,w_q)\}_{q=1}^{N_{\mathrm{ip},f}}$ to obtain
\begin{equation}
  \label{eq:rhs_quadrature}
  b_i \approx \frac12\sum_{q=1}^{N_{\mathrm{ip},f}} w_q\, g^f_{L,q}\,N_i\!\big(x_L(r_q)\big)
          + \frac12\sum_{q=1}^{N_{\mathrm{ip},f}} w_q\, g^f_{R,q}\,N_i\!\big(x_R(r_q)\big),
\end{equation}
where $g^f_{L,q}$ and $g^f_{R,q}$ denote the fine values at the $q$-th Gauss point on the left and right child, respectively.

Finally, the mass matrix $M_{ij}$ can be evaluated as:
\begin{equation}
  \label{eq:diag_mass}
  M_{ij} \approx \sum_{k=1}^{N_{\mathrm{ip},c}} w_k\,N_i(\vec{x}_{k}^c)N_j(\vec{x}_{k}^c) = w_i\,\delta_{ij},
\end{equation}
We note that the choices of the basis functions with the Lagrange property at the quadrature points lead to a diagonal mass matrix. Combining \eqnref{eq:rhs_quadrature} and \eqnref{eq:diag_mass}, we can directly compute the coarse coefficients without solving a linear system.

\textbf{Restriction matrix form for 1D:}
The above procedure can be compactly expressed in matrix--vector form. Collect the fine values on the two children into a single vector
\[
g^f :=
\big[g^f_{L,1},\ldots,g^f_{L,N_{\mathrm{ip},f}},\; g^f_{R,1},\ldots,g^f_{R,N_{\mathrm{ip},f}}\big]^T \in \mathbb{R}^{2\mathrm{N}_{\mathrm{ip},f}}.
\]
Then the projection restriction can be written as a matrix--vector product
\begin{equation}
\label{eq:restriction_matrix_form}
g_c = R^{1D}\, g^f,
\qquad
g_c = [g_{c_1},\ldots,g_{c_{N_{\mathrm{ip},c}}}]^T \in \mathbb{R}^{N_{\mathrm{ip},c}},
\end{equation}
where the restriction matrix $R^{1D}\in\mathbb{R}^{N_{\mathrm{ip},c}\times 2N_{\mathrm{ip},f}}$ has entries
\begin{equation}
\label{eq:R_entries}
R^{1D}_{i,\,(c,q)}
=
\frac{1}{w_i}\left(\frac12\, w_q\right)\,
N_i\!\big(x_c(r_q)\big),
\qquad
c\in\{L,R\},\;\; q=1,\ldots,N_{\mathrm{ip},f}.
\end{equation}
Here $x_L(r)=\tfrac12(r-1)$ and $x_R(r)=\tfrac12(r+1)$ are the child-to-parent maps, and the factor $\tfrac12$ is the Jacobian $dx=\tfrac12 dr$.

\subsection{Extending to arbitrary order of integration}
Certain applications require a larger number of integration points ($\geq (p+1)$) for better accuracy (e.g., \cite{saurabh2021industrial}). The above construction can be generalized to arbitrary integration order and with different numbers of quadrature points on fine and coarse elements. The local $L^2$ projection framework remains unchanged: we still compute the coarse quadrature values by solving the local mass-matrix system \eqnref{eq:mass_system}, where the right-hand side $b_i$ is computed by summing contributions from all children using appropriate quadrature rules. The only difference is that the mass matrix $M_{ij}$ is no longer diagonal as in \eqnref{eq:rhs_quadrature}, and we must solve the linear system to obtain the coarse quadrature values. We also note that these can be precomputed for a given choice of basis and quadrature rules and do not depend on the actual fine-mesh data, but rather on the position and weights of the quadrature points. The restriction matrix form \eqnref{eq:restriction_matrix_form}--\eqnref{eq:R_entries} also remains valid, with the entries modified to account for the different numbers of quadrature points on fine and coarse elements. The remainder of this work focuses on the use of standard Gauss--Legendre quadrature with $(p+1)$ points for $Q_p$ elements.

\subsection{Extending to multiple dimensions}
The 1D restriction matrix $R^{1D}$ defined in \eqnref{eq:R_entries} can be extended to multiple dimensions via tensor products. For example, in 2D, the restriction matrix is given by
\begin{equation}
\label{eq:restriction_matrix_2D}
R^{2D} = R^{1D} \otimes R^{1D},
\end{equation}
where $\otimes$ denotes the Kronecker product. The entries of $R^{2D}$ can be explicitly written as
\begin{equation}
\label{eq:R_entries_2D}
R^{2D}_{(i_1,i_2),\,(c_1,q_1,c_2,q_2)}
=
R^{1D}_{i_1,\,(c_1,q_1)}\,
R^{1D}_{i_2,\,(c_2,q_2)},
\end{equation}
where $i_1,i_2=1,\ldots,N_{\mathrm{ip},c}$, $c_1,c_2\in\{L,R\}$, and $q_1,q_2=1,\ldots,N_{\mathrm{ip},f}$. The multi-dimensional restriction is then given by
\begin{equation}
\label{eq:restriction_matrix_form_2D}
g_c = R^{2D}\, g^f,
\end{equation}
where $g^f$ collects the fine values on all children in lexicographic order. We note that we do not need to explicitly form the multi-dimensional restriction matrix; instead, we can apply the Kronecker structure on-the-fly during the restriction operation.

\subsection{Recovering coarse nodal DOFs via global $L^2$ projection}
\label{subsec:l2_projection}

The values $g_c(\vec{x}_j^c)$ obtained above live at quadrature points on the coarse element. To obtain the coarse-mesh nodal degrees of freedom, we compute the $L^2$ projection of $g_c$ onto the coarse CG finite element space. To assemble the mass matrix and right-hand side in \eqnref{eqn:mass_matrix_system}, the integration is performed using the same quadrature rule that was used to define the coarse quadrature-point values $g_c$ (\eqnref{eq:R_entries}).

Let $V_h \subset L^2(\Omega)$ denote the coarse CG finite element space with Lagrange basis functions $\{N_a\}_{a=1}^{n_h}$ (of degree $p\ge 1$). We seek $g^h \in V_h$ that minimizes the $L^2$ error with respect to the quadrature-point field $g_c$:
\begin{equation}
    g^h
    \;=\;
    \arg\min_{v^h \in V_h}
    \frac{1}{2}\,\|v^h - g_c\|_{L^2(\Omega)}^2.
    \label{eqn:l2_min}
\end{equation}
The associated variational condition is: find $g^h \in V_h$ such that
\begin{equation}
    (g^h, w^h) = (g_c, w^h)
    \qquad \forall\, w^h \in V_h,
    \label{eqn:l2_min_var}
\end{equation}
where $(\cdot,\cdot)$ denotes the standard $L^2$ inner product. Writing $g^h(\vec{x})=\sum_{a=1}^{n_h} G_a\,N_a(\vec{x})$ yields the familiar mass-matrix system
\begin{equation}
    \mathbf{M}\,\mathbf{G} = \mathbf{b},
    \qquad
    M_{ab} = (N_a,N_b),
    \qquad
    b_a = (g_c,N_a),
    \label{eqn:mass_matrix_system}
\end{equation}
which is solved to obtain the nodal coefficients $\mathbf{G}=\{G_a\}$ on the coarse mesh. In practice, the right-hand side $(g_c,N_a)$ is assembled using the coarse quadrature points and the values $g_c(\vec{x}_j^c)$. We note that while the dimension of the mass matrix $\mathbf{M}$ depends only on the number of coarse-mesh degrees of freedom, the integrals for assembling both $\mathbf{M}$ and $\mathbf{b}$ must be evaluated using the same quadrature rule that was used to obtain the coarse-point values $g_c$.

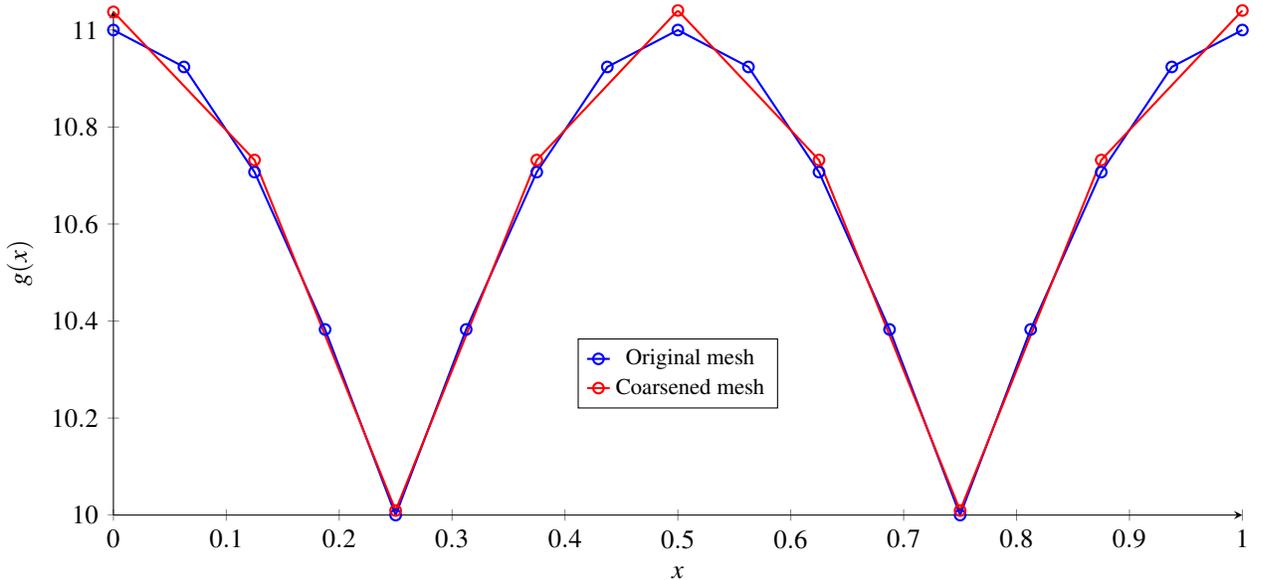
\begin{figure}[h]
\begin{subfigure}[b]{1.0\linewidth}
	\centering
\begin{tikzpicture}
\begin{axis}[
          width=1.0\linewidth, 
          height = 0.5\linewidth,
            xlabel={$x$},
    ylabel={$g(x)$},
    axis x line=bottom,
    axis y line=left,
	legend style={at={(0.50,0.35)},anchor= north,nodes={scale=0.95, transform shape}},
]
     \addplot[mark=o,color = blue,  thick,name path=A] 
        table[x expr={\thisrow{x}))},y expr=(\thisrow{val},col sep=space]{Data/Curve1D/fineMesh.txt};
     \addplot[mark=o,color = red,  thick,name path=B] 
        table[x expr={\thisrow{x}))},y expr=(\thisrow{val},col sep=space]{Data/Curve1D/Projection.txt};
        \legend{\small{Original mesh}, \small{Coarsened mesh}}
  \end{axis}
  \end{tikzpicture}
  \end{subfigure}
  \caption{\textbf{Mass conserving interpolation:} Figure showing the redistribution of $\phi$ on the coarse mesh when the solution is projected from the fine mesh, such that $\int_{\Omega} g(x) d\Omega$ remains constant. 	$\bigg[\int_{\Omega} g(x) \; d\Omega \bigg]_{\mathcal{M}_O}$ =  
  	$\bigg[\int_{\Omega} g(x) \; d\Omega \bigg]_{\mathcal{M}_C}$ = 10.6284 , where $\mathcal{M}_O$ is the original mesh and $\mathcal{M}_C$ is the coarsened mesh.}
  	\label{fig:projection}
\end{figure}

\begin{figure}[h]
\begin{subfigure}[b]{0.5\linewidth}
	\centering
\begin{tikzpicture}
\begin{axis}[
          width=\linewidth, 
          height = \linewidth,
            xlabel={$x$},
    ylabel={$g(x)$},
    axis x line=bottom,
    axis y line=left,
	legend style={at={(0.50,0.30)},anchor= north,nodes={scale=0.95, transform shape}},
]
     \addplot[mark=none,color = blue,   thick] 
        table[x expr={\thisrow{Points_0}},y expr=\thisrow{ch},col sep=comma]{Data/1D_Quad/refine.csv};
    \addplot[mark=none,color = red,  thick] 
        table[x expr={\thisrow{Points_0}},y expr=\thisrow{ch},col sep=comma]{Data/1D_Quad/injects.csv};
 \addplot[only marks, mark=o, color=blue, ultra thick]
        table[x expr={\thisrow{Points_0}},y expr=\thisrow{ch},col sep=comma]{Data/1D_Quad/points_refine.csv};

 \addplot[only marks, mark=square, color=red, ultra thick]
        table[x expr={\thisrow{Points_0}},y expr=\thisrow{ch},col sep=comma]{Data/1D_Quad/points_inject.csv};
\legend{Original mesh, Coarse mesh}
  \end{axis}
  \end{tikzpicture}
  \caption{Coarsening by injection, $\bigg[\int_\Omega g(x) dx\bigg]_{\mathcal{M}_C}$ = 10.6381}
  \end{subfigure}
  \begin{subfigure}[b]{0.5\linewidth}
	\centering
\begin{tikzpicture}
\begin{axis}[
          width=\linewidth, 
          height = \linewidth,
            xlabel={$x$},
    ylabel={$g(x)$},
    axis x line=bottom,
    axis y line=left,
	legend style={at={(0.50,0.30)},anchor= north,nodes={scale=0.95, transform shape}},
]
     \addplot[mark=none,color = blue,  thick] 
        table[x expr={\thisrow{Points_0}},y expr=\thisrow{ch},col sep=comma]{Data/1D_Quad/refine.csv};
        \addplot[mark=none,color = cpu1,   thick] 
        table[x expr={\thisrow{Points_0}},y expr=\thisrow{ch},col sep=comma]{Data/1D_Quad/L2.csv};
 \addplot[only marks, mark=o, color=blue, ultra thick]
        table[x expr={\thisrow{Points_0}},y expr=\thisrow{ch},col sep=comma]{Data/1D_Quad/points_refine.csv};

 \addplot[only marks, mark=x, color=cpu1, ultra thick, mark size=4pt]
        table[x expr={\thisrow{Points_0}} ,y expr=\thisrow{ch},col sep=comma]{Data/1D_Quad/points_L2.csv};
\legend{Original mesh, Coarse mesh}
  \end{axis}
  \end{tikzpicture}
\caption{Coarsening by mass-conserving interpolation, $\bigg[\int_\Omega g(x) dx\bigg]_{\mathcal{M}_C}$ = 10.6367}
  \end{subfigure}
  \caption{\textbf{Interpolation using quadratic basis function}: Figure comparing the coarsening of a function $g(x) = |\cos(2\pi x)| + 10 $ defined on a fine mesh with 8 quadratic elements to a coarse mesh with 4 quadratic elements using (a) injection and (b) mass-conserving interpolation. The integral of the function over the domain $\Omega=[0,1]$ on the original refined mesh (=10.6367) is preserved up to numerical precision using the proposed mass-conserving interpolation scheme, whereas injection results in a noticeable deviation.}
  \label{fig:Coarsening1D_quad}
\end{figure}
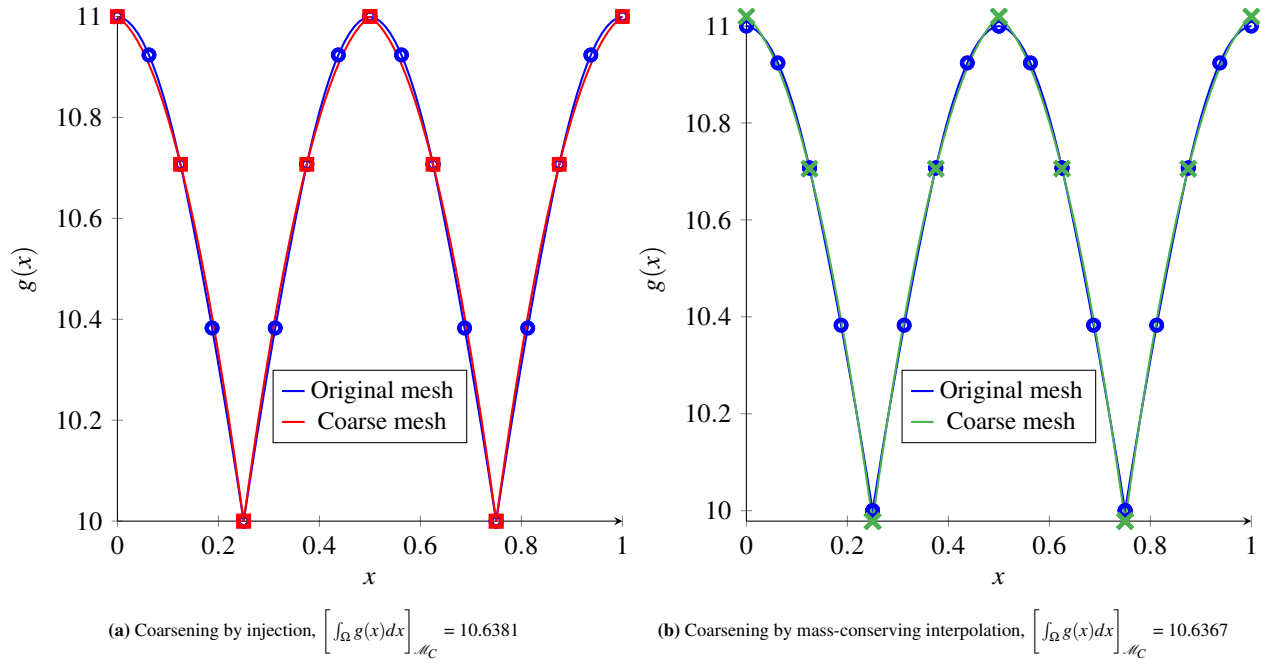

\textbf{Interpretation:} The overall procedure can be interpreted as a \emph{global redistribution} of the conserved quantity: coarse nodal values may differ locally from those produced by injection, but the integral constraint \eqnref{eqn:mass_consv_disc} is satisfied by construction. This is expected as the basis functions representing the quantities are global and continuous.  \figref{fig:projection} illustrates this effect for $g(x)=|\cos(2\pi x)|+10$ transferred from a fine mesh to a coarse mesh: local nodal values differ, while the global integral is preserved to numerical precision. Similarly. \figref{fig:Coarsening1D_quad} compares the standard injection-based coarsening with the proposed conservative coarsening for the quadratic basis function. While injection-based coarsening introduces noticeable mass drift, the proposed conservative scheme maintains mass conservation to numerical precision.

\section{Algorithms}
\label{sec:Algorithm}

In this section, we describe the algorithmic steps required to perform mass-conserving interpolation during the coarsening stage of AMR based on the numerical procedure outlined in \secref{sec:Numerics}. We recall that the refinement and coarsening stage is mutually exclusive and we focus only on the coarsening stage here, i.e., elements marked now only contain \texttt{COARSEN} or \texttt{NO\_CHANGE} flags. To enable coarsening of elements, all children of a particular parent must meet the following two criteria:

\begin{enumerate}
\item All child must be at the same level.
\item All child must have \texttt{COARSEN} flags.
\end{enumerate}

\def \cs{0.07}
\tikzexternaldisable
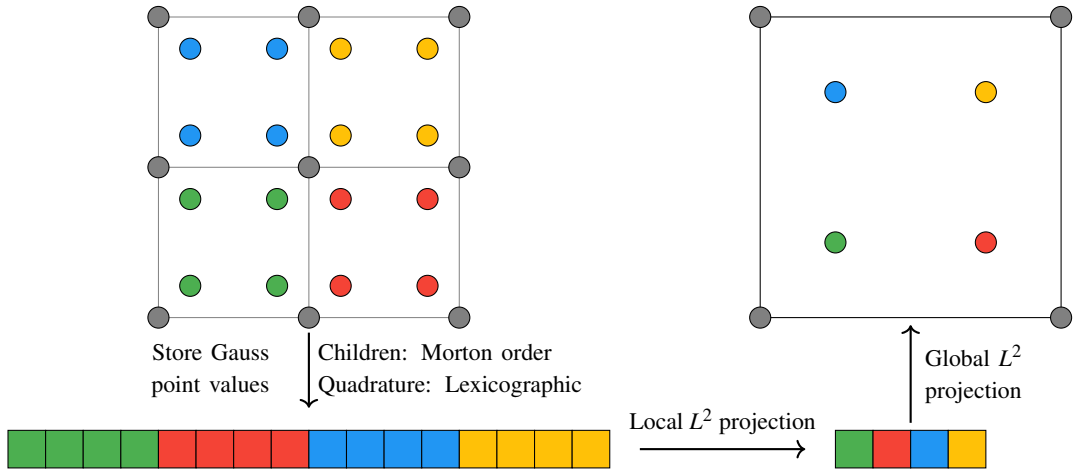
\begin{figure}
\centering
    \begin{tikzpicture}[scale=2]
    \draw[step=1.0cm,gray,very thin] (-3.0,-2) grid (-1.0,0);
    
    \draw[fill=cpu1](-3 + 1*0.2113248654,-2.0 + 1*0.2113248654) circle (\cs cm);
    \draw[fill=cpu1](-3 + 1*0.2113248654,-2.0 + 1*0.78867513459) circle (\cs cm);
    \draw[fill=cpu1](-3 + 1*0.78867513459,-2.0 + 1*0.2113248654) circle (\cs cm);
    \draw[fill=cpu1](-3 + 1*0.78867513459,-2.0 + 1*0.78867513459) circle (\cs cm);

    \draw[fill=cpu2](-2 + 1*0.2113248654,-1.0 + 1*0.2113248654) circle (\cs cm);
    \draw[fill=cpu2](-2 + 1*0.2113248654,-1.0 + 1*0.78867513459) circle (\cs cm);
    \draw[fill=cpu2](-2 + 1*0.78867513459,-1.0 + 1*0.2113248654) circle (\cs cm);
    \draw[fill=cpu2](-2 + 1*0.78867513459,-1.0 + 1*0.78867513459) circle (\cs cm);

    \draw[fill=cpu3](-2 + 1*0.2113248654,-2.0 + 1*0.2113248654) circle (\cs cm);
    \draw[fill=cpu3](-2 + 1*0.2113248654,-2.0 + 1*0.78867513459) circle (\cs cm);
    \draw[fill=cpu3](-2 + 1*0.78867513459,-2.0 + 1*0.2113248654) circle (\cs cm);
    \draw[fill=cpu3](-2 + 1*0.78867513459,-2.0 + 1*0.78867513459) circle (\cs cm);

    \draw[fill=cpu4](-3 + 1*0.2113248654,-1.0 + 1*0.2113248654) circle (\cs cm);
    \draw[fill=cpu4](-3 + 1*0.2113248654,-1.0 + 1*0.78867513459) circle (\cs cm);
    \draw[fill=cpu4](-3 + 1*0.78867513459,-1.0 + 1*0.2113248654) circle (\cs cm);
    \draw[fill=cpu4](-3 + 1*0.78867513459,-1.0 + 1*0.78867513459) circle (\cs cm);

	\draw[fill=gray](-3 + 1*0.0,-2.0 + 1*0.0) circle (\cs cm);
	\draw[fill=gray](-3 + 1*0.0,-1.0 + 1*0.0) circle (\cs cm);
	\draw[fill=gray](-3 + 1*0.0,-0.0 + 1*0.0) circle (\cs cm);
	
	\draw[fill=gray](-2 + 1*0.0,-2.0 + 1*0.0) circle (\cs cm);
	\draw[fill=gray](-2 + 1*0.0,-1.0 + 1*0.0) circle (\cs cm);
	\draw[fill=gray](-2 + 1*0.0,-0.0 + 1*0.0) circle (\cs cm);

	\draw[fill=gray](-1 + 1*0.0,-2.0 + 1*0.0) circle (\cs cm);
	\draw[fill=gray](-1 + 1*0.0,-1.0 + 1*0.0) circle (\cs cm);
	\draw[fill=gray](-1 + 1*0.0,-0.0 + 1*0.0) circle (\cs cm);	

    \draw (1,-2) rectangle (3,0);
    
    \draw[fill=cpu1](1 + 1*0.5,-2.0 + 1*0.5) circle (\cs cm);
    \draw[fill=cpu4](1 + 1*0.5,-1.0 + 1*0.5) circle (\cs cm);
    \draw[fill=cpu3](2 + 1*0.5,-2.0 + 1*0.5) circle (\cs cm);
    \draw[fill=cpu2](2 + 1*0.5,-1.0 + 1*0.5) circle (\cs cm);
    \draw[->, thick] (-2,-2.1) -- (-2,-2.6);
    \node[text width=2cm, align=right, anchor=east] at (-2.2,-2.35) {
        \small{Store Gauss\\point values}
    };
    
    \node[text width=4cm, align=left, anchor=west] at (-2.0,-2.35) {
        \small{Children: Morton order\\Quadrature: Lexicographic}
    };
    \foreach \x in {0,0.25,0.5,0.75}
        \draw[fill=cpu1] (-4 + \x,-3) rectangle (-3.75+\x,-2.75);

    \foreach \x in {0,0.25,0.5,0.75}
        \draw[fill=cpu3] (-3 + \x,-3) rectangle (-2.75+\x,-2.75);
    
    \foreach \x in {0,0.25,0.5,0.75}
        \draw[fill=cpu4] (-2 + \x,-3) rectangle (-1.75+\x,-2.75);
        
    \foreach \x in {0,0.25,0.5,0.75}
        \draw[fill=cpu2] (-1 + \x,-3) rectangle (-0.75+\x,-2.75);

    \draw[->, thick] (0.2,-2.875) -- (1.3,-2.875);
    \node[above=2pt] at (0.75,-2.875) {\small Local $L^2$ projection};
    \draw[fill=cpu1] (1.5,-3) rectangle (1.75,-2.75);
    \draw[fill=cpu3] (1.75,-3) rectangle (2.0,-2.75);
    \draw[fill=cpu4] (2.0,-3) rectangle (2.25,-2.75);
    \draw[fill=cpu2] (2.25,-3) rectangle (2.50,-2.75);
    
    \draw[fill=gray](1 + 1*0.0,-2.0 + 0*0.0) circle (\cs cm);
    \draw[fill=gray](1 + 1*0.0,-0.0 + 0*0.0) circle (\cs cm);
    
    \draw[fill=gray](3 + 1*0.0,-2.0 + 0*0.0) circle (\cs cm);	
    \draw[fill=gray](3 + 1*0.0,-0.0 + 0*0.0) circle (\cs cm);

    \draw[->, thick] (2.0,-2.7) -- (2.0,-2.05);
    \node[right=2pt, text width=2cm, align=left] at (2.0,-2.375) {\small Global $L^2$ projection};
    
    \end{tikzpicture}
    \caption{\textit{Interpolation scheme:} The elements of the fine mesh  (left figure) is coarsened one time  to result in a single element shown on the right-hand side. We first store the values at all the Gauss points (4 in 2D and 8 in 3D) in a single vector. The children are arranged in Morton order, whereas the integration points within each elements are arranged in lexicographic order. The values at these Gauss points are then interpolated using local $L^2$ projection described in \Algref{alg:l2_restrict_dd}. Finally, we perform an $L^2$ projection from the Gauss points to the nodal points to get the final values at the nodes of the coarse mesh.}
    \label{fig:interpolation}
\end{figure}
\tikzexternalenable

\begin{algorithm}[t]
\caption{Local $L^2$ restriction from fine to coarse grid on one coarsened element}
\label{alg:l2_restrict_dd}
\begin{algorithmic}[1]
\Require Dimension $dim\in\{2,3\}$,  $n$ number of integration point per dimension
\Require 1D restriction matrix $R^{1D}$ 
\Require Gauss point values on one child $c$ at level $l{+}1$: $g^f[\texttt{f\_itg}]$, $\texttt{f\_itg}=0,\ldots,N_{ip,f}-1$
\Require Child id: $\texttt{child}\in\{0,\ldots,2^{dim}-1\}$  \Comment {Child index in Morton order}
\Ensure Coarse samples on parent: $U[\texttt{c\_dof}]$, $\texttt{c\_dof}=0,\ldots,N_{ip,c}-1$

\State  $(c_x,c_y,c_z) \gets$ Decode\_Morton(\texttt{child},dim) \Comment{\Algref{alg:decode_morton}}
\State  $U\gets 0$ \Comment{Initialize coarse Gauss point values to zero}

\For{\texttt{c\_itg} $=0$ to $N_{ip,c}-1$}  \Comment{coarse point Gauss point values in lexicographic order}
    \State $(I_x,I_y,I_z) \gets$ Decode\_Lex(\texttt{c\_itg},dim,n) \Comment{\Algref{alg:decode_lex}}
     \State $\texttt{acc} \gets 0$
    \For{\texttt{f\_itg} $=0$ to $N_{ip,f}-1$} \Comment{fine point Gauss point in lexicographic order}
        \State  $\mathbf{\ell}=(\ell_x,\ell_y,\ell_z) \gets$  Decode\_Lex(\texttt{f\_itg},dim,n) \Comment{\Algref{alg:decode_lex}}
        \State w $\gets 1$
        \For{d $\in$ \{x,y,z\}}
            \State $i_d \gets c_d\,n + \ell_d$ 
            \State w $\gets$ w $\cdot  R^{1D}\texttt{[I\_d][i\_d]}$  \Comment{Kronecker-product}
        \EndFor
        \State $\texttt{acc} \gets \texttt{acc} + w \cdot g^f[\texttt{f\_itg}]$
    \EndFor
    \State $U[\texttt{c\_itg}] \gets U[\texttt{c\_itg}] + \texttt{acc}$
\EndFor
\State \Return $U$
\end{algorithmic}
\end{algorithm}

\begin{algorithm}[t]
\caption{Decode\_lex :  Decode lexicographic index to multi-dimensional coordinates}
\label{alg:decode_lex}
\begin{algorithmic}[1]
\Require Lexicographic index $\texttt{lex\_idx}$
\Require Dimension $dim\in\{2,3\}$
\Require Number of points per dimension: $n$
\Ensure Coordinates $(I_x, I_y, I_z)$

\If{$dim = 2$}
    \State $I_x \gets \texttt{lex\_idx} \mod n$
    \State $I_y \gets \lfloor \texttt{lex\_idx} / n \rfloor$
    \State $I_z \gets 0$
\ElsIf{$dim = 3$}
    \State $I_x \gets \texttt{lex\_idx} \mod n$
    \State $I_y \gets \lfloor (\texttt{lex\_idx} / n) \rfloor \mod n$
    \State $I_z \gets \lfloor \texttt{lex\_idx} / n^2 \rfloor$
\EndIf
\State \Return $(I_x, I_y, I_z)$
\end{algorithmic}
\end{algorithm}

\begin{algorithm}[t]
\caption{Decode\_Morton : Decode Morton (Z-order) index to child coordinates}
\label{alg:decode_morton}
\begin{algorithmic}[1]
\Require Child index $\texttt{child}\in\{0,\ldots,2^{dim}-1\}$ in Morton order
\Require Dimension $dim\in\{2,3\}$
\Ensure Child coordinates $(c_x, c_y, c_z)\in\{0,1\}^{dim}$

\State $c_x \gets \texttt{child} \,\&\, 1$ \Comment{Extract bit 0}
\State $c_y \gets (\texttt{child} \,\&\, 2) \gg 1$ \Comment{Extract bit 1}
\If{$dim = 3$}
    \State $c_z \gets (\texttt{child} \,\&\, 4) \gg 2$ \Comment{Extract bit 2}
\Else
    \State $c_z \gets 0$
\EndIf
\State \Return $(c_x, c_y, c_z)$
\end{algorithmic}
\end{algorithm}
\figref{fig:interpolation} briefs the step required for coarsening. The main steps for elements that are classified as coarsening can be stated as follows:
\begin{enumerate}
	\item \textbf{Gauss point interpolation:} The values at the Gauss points are interpolated from the nodal values for each element using the associated basis functions. 
	\item  \textbf{Transfer to the coarser element:} Once we have the values at the Gauss points for all the children elements, we need to transfer these values to the coarser parent element. This is achieved by performing a local $L^2$ projection from the fine mesh to the coarse mesh as described in \Algref{alg:l2_restrict_dd}. First, we store the values at all the Gauss points of the children elements in a single vector. The children are arranged in Morton order, whereas the integration points within each element are arranged in lexicographic order. Next, we perform the local $L^2$ projection from the fine mesh to the coarse mesh to evaluate the values at the Gauss points of the coarser element. \figref{fig:interpolation} illustrates this process in 2D.

    \item \textbf{$L^2$ Projection:} After the previous step, we have the values at the Gauss points for the coarser element. To determine the values at the nodal points, we perform global $L^2$ projection from the Gauss points to the nodal points, as described in \secref{sec:Numerics}.

\end{enumerate}

Next, we provide a detailed explanation of the algorithm and data structure required to perform mass-conserving interpolation. We begin by allocating two vectors of size \texttt{elem\_size} $\times$ \texttt{ngp} for the mesh before ($\mathcal{M}_1$) and after the AMR step ($\mathcal{M}_2$). Here, \texttt{elem\_size} represents the number of elements and \texttt{ngp} denotes the number of Gauss points used to evaluate the field variable. We assume (without loss of generalization) that the same number of Gauss points is employed per element. Basically, $N_{ip,c}  = N_{ip,f}$ =  \texttt{ngp} = $(p+1)^{dim}$.

We iterate over each Gauss point in the original mesh $\mathcal{M}_1$ and store its corresponding value in the vector. For elements that do not undergo the coarsening stage, the values associated with those elements are copied into the vector corresponding to the mesh $\mathcal{M}_2$. However, if the cells do undergo coarsening, we compute the local $L^2$ projection from the fine mesh to the coarse mesh as described in \Algref{alg:l2_restrict_dd} to transfer the values from the $\mathcal{M}_1$ vector. This operation is completely local and independent for each coarsened element, involving only the children elements and their corresponding entries in the fine elements. The results of this local projection are then stored in the appropriate locations within the coarse element vector. This step can be thought of as transferring the values at the cell level (cell intergrid transfer). Once the values at the Gauss points for the coarser mesh $\mathcal{M}_2$ are obtained, we proceed to perform an $L^2$ projection from the Gauss points to the nodal points by inverting the mass matrix, as described in \secref{sec:Numerics}.

\section{Results}
\label{sec:results}

In this section, we compare the proposed \emph{field-conserving octree coarsening} operator against the standard injection-based child-to-parent transfer. We first perform manufactured solution tests to verify convergence and conservation properties in \secref{subsec:results_mms}. Next, we evaluate the schemes in the context of phase separation phenomena governed by the Cahn--Hilliard and Cahn--Hilliard--Navier--Stokes equations in \secref{subsec:results_ch} and \secref{sec:bubble_case2}, respectively.
We study spinodal decomposition governed by the Cahn--Hilliard (CH) equation using both polynomial and Flory--Huggins free energies, and then extend the evaluation to a representative Cahn--Hilliard--Navier--Stokes (CHNS) multiphase flow test case. Throughout, the mesh is adapted with \emph{single-level refinement/coarsening on balanced octrees}, so that differences between runs can be attributed primarily to the coarsening transfer operator.

\subsection{Convergence Error: Method of Manufactured Solution}
\label{subsec:results_mms}
To verify the accuracy of the proposed conservative coarsening operator, we first perform a convergence study using the method of manufactured solutions (MMS). We consider the diffusion equation:
\begin{equation}
    \frac{\partial \phi}{\partial t} = \nabla \cdot (\kappa \nabla \phi)
    \label{eq:mms_diffusion}
\end{equation}
on domain $\Omega = [0,1]^2$, with homogeneous Neumann boundary conditions:
\begin{equation}
    \nabla \phi \cdot \hat{n} = 0 \quad \text{on } \partial \Omega,
\end{equation}
where $\kappa$ is the diffusion coefficient. We chose the manufactured solution as:
\begin{equation}
    \phi(x,y,t) = 1 + A \cos(2 \pi x) \cos(2 \pi  y) e^{-\kappa 8 \pi^2 t}
    \label{eq:mms_solution}
\end{equation}
with parameters $A=0.1$ and $\kappa=0.03$. The solution of the above PDE is mass-conserving, i.e., $\int_\Omega \phi d(\Omega)$ remains constant over time and equals the initial mass $\int_\Omega \phi(t=0) d\Omega = 1$. 

To numerically simulate the PDE, we use second-order Crank-Nicholson time integration with timestep $\Delta t = 0.01$ for linear basis and $\Delta t = 0.001$ for quadratic basis function. We chose a lower timestep for the quadratic basis function to ensure that the temporal error is negligible compared to the spatial error. The final time is set to $t_f = 1.0$.

To perform the convergence study, we consider a sequence of uniformly refined meshes at different octree levels. At each step, we compute the local gradient of the solution $\nabla \phi$ and adapt the mesh based on the following criterion:
\begin{equation}
    \eta_e = \int_{\Omega_e} \|\nabla \phi \|_{L_2(e)} \mathrm{d}\Omega_e< \tau,
    \label{eq:mms_adapt_criteria}
\end{equation}
where $\tau$ is a user-defined mesh-dependent parameter that scales with mesh size h. \tabref{tab:mass_drift_mms} reports the value for different level. We choose $\eta_e$ as a simple heuristic based on the derivatives evaluated at quadrature point. Elements satisfying \eqnref{eq:mms_adapt_criteria} are coarsened, whereas the others are marked as \texttt{NO\_CHANGE}. Additionally, at each iteration, we coarsen a fixed percentage (10\%) of the total elements by sorting them according to \eqnref{eq:mms_adapt_criteria} and coarsening the elements with the lowest values. In each simulation, we consider the octree elements to be at level $l$ or $l-1$. An octant at level $l$ is coarsened to level $l-1$ if it satisfies the above criteria, while the octant already at level $l-1$ is unchanged.
 This ensures that we have sufficient coarsening operations to evaluate the proposed scheme. We perform this adaptation procedure at every timestep throughout the simulation.

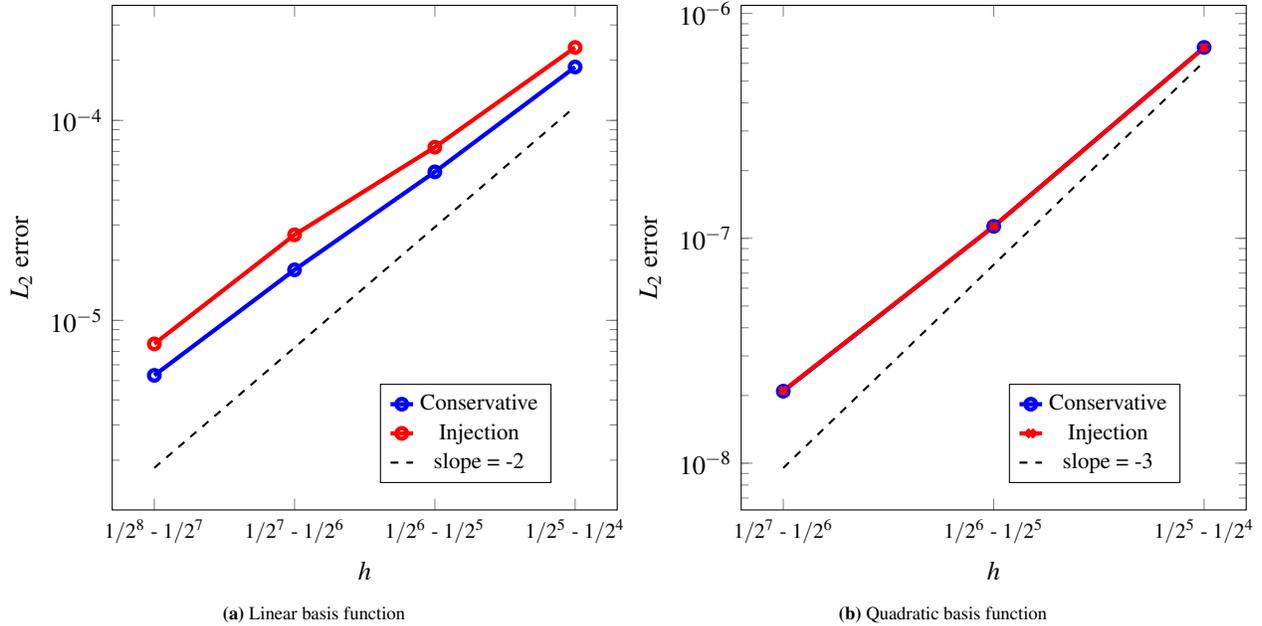
\begin{figure}
    \begin{subfigure}[b]{0.5\linewidth}
	\centering
\begin{tikzpicture}
\begin{loglogaxis}[
          width=1\linewidth, 
          height = 1\linewidth,
            xlabel={$h$},
    ylabel={$L_2$ error},
          xtick={1/2^5, 1/2^6, 1/2^7,1/2^8},
         xticklabels={\footnotesize{$1/2^5$ - $1/2^4$}, \footnotesize{$1/2^6$ - $1/2^5$}, \footnotesize{$1/2^7$ - $1/2^6$},\footnotesize{$1/2^8$ - $1/2^7$}},
	legend style={at={(0.70,0.25)},anchor= north,nodes={scale=0.95, transform shape}},
]
        \addplot[mark=o, color=blue, ultra thick] table[
            x expr={1/2^(\thisrow{Level})}, 
            y expr={\thisrow{L2_Error}}, 
            col sep=space
        ] {Data/MMS/Error_linear.txt};

        \addplot[mark=o, color=red, ultra thick] table[
            x expr={1/2^(\thisrow{Level})}, 
            y expr={\thisrow{Inj_Error}}, 
            col sep=space
        ] {Data/MMS/Error_linear.txt};
\addplot[black, dashed, thick, domain=1/256:1/32] {0.12*x^2};        
     \legend{\small{Conservative }, \small {Injection} ,\small{slope = -2}}
  \end{loglogaxis}
  \end{tikzpicture}
  \caption{Linear basis function}
  \end{subfigure}
  \begin{subfigure}[b]{0.5\linewidth}
	\centering
\begin{tikzpicture}
\begin{axis}[
          width=1\linewidth, 
          height = 1\linewidth,
            xlabel={$h$},
    ylabel={$L_2$ error},
    xmode=log,
    ymode=log,
     xtick={1/2^5, 1/2^6, 1/2^7},
          xticklabels={\footnotesize{$1/2^5$ - $1/2^4$}, \footnotesize{$1/2^6$ - $1/2^5$}, \footnotesize{$1/2^7$ - $1/2^6$}},
    legend style={at={(0.70,0.25)},anchor= north,nodes={scale=0.95, transform shape}},
]
        \addplot[mark=o, color=blue, ultra thick] table[
            x expr={1/2^(\thisrow{Level})}, 
            y expr={\thisrow{L2_Error}}, 
            col sep=space
        ] {Data/MMS/Error_quad.txt};

        \addplot[mark=x, color=red, ultra thick] table[
            x expr={1/2^(\thisrow{Level})}, 
            y expr={\thisrow{Inj_Error}}, 
            col sep=space
        ] {Data/MMS/Error_quad.txt};
       \addplot[black, dashed, thick, domain=1/128:1/32] {0.02*x^3};
           \legend{\small{Conservative }, \small{Injection},\small{slope = -3}}
  \end{axis}
  \end{tikzpicture}
  \caption{Quadratic basis function}
  \end{subfigure}
  \caption{\textbf{MMS error:} Figure comparing the $L_2$ error for linear and quadratic basis function. Each run contains a mesh at two different resolutions, differing by one level, as shown by the range in $h$. Both coarsening schemes exhibit approximately optimal convergence rates, with the proposed mass-conserving scheme showing slightly lower errors for linear basis functions and comparable errors for quadratic basis functions.}
    \label{fig:MMS_convergence}
\end{figure}

\figref{fig:MMS_convergence} shows  $L_2$ error in $\phi$ at final time $t_f=1.0$ for linear and quadratic basis functions. We do not observe a significant difference in convergence rate between the two interpolation schemes. Both schemes exhibit optimal convergence rates of approximately 2 for linear basis functions and approximately 3 for quadratic basis functions, as expected. We observe a slight difference in error magnitude between the two schemes, with the proposed conservative coarsening operator yielding marginally lower errors for linear basis functions and almost identical errors for quadratic basis functions. This demonstrates that the proposed conservative coarsening operator maintains the expected accuracy of the numerical solution without introducing additional errors during mesh adaptation.

\begin{table}[h]
\centering
\begin{tabular}{cc|cc|cc}
\hline
\multirow{2}{*}{Level} & \multirow{2}{*}{$\tau$} & \multicolumn{2}{c|}{Linear} & \multicolumn{2}{c}{Quadratic} \\
\cline{3-6}
 & & Injection & Conservative & Injection & Conservative \\
\hline
5 & 1 $\times 10^{-2}$
 & $3.72 \times 10^{-7}$ & $2.08 \times 10^{-13}$ & $1.71 \times 10^{-10}$ & $3.60 \times 10^{-13}$ \\
6 & 5 $\times 10^{-3}$ & $1.97 \times 10^{-9}$ & $8.84 \times 10^{-13}$ & $9.75 \times 10^{-12}$ & $1.16 \times 10^{-12}$ \\
7 & 2.5 $\times 10^{-3}$ & $6.39 \times 10^{-10}$ & $3.83 \times 10^{-13}$ & $4.00 \times 10^{-12}$ & $9.75 \times 10^{-13}$ \\
8 & 1.25 $\times 10^{-3}$ & $1.79 \times 10^{-10}$ & $2.75 \times 10^{-13}$ & -- & -- \\
\hline
\end{tabular}
\caption{Mass drift for the MMS test case at final time $t_f=1.0$ for different octree levels using the injection and conservative coarsening scheme.}
\label{tab:mass_drift_mms}
\end{table}
We finally compute the mass drift in $\phi$ at final time $t_f=1.0$. The mass drift is defined as:
\begin{equation*}
    \Delta \phi = \left| \int_{\Omega} \phi (t = t_f) d\Omega - \int_{\Omega} \phi (t = 0) d\Omega \right|
\end{equation*} 
This mass drift can be attributed to the coarsening operations performed during mesh adaptation. On the uniform mesh without any adaptation, the numerical solution conserves $\int \phi$ up to machine precision. However, when coarsening is performed, the injection-based scheme does not guarantee mass conservation, leading to a non-zero mass drift. \tabref{tab:mass_drift_mms} reports the mass drift at final time $t_f=1.0$ for different octree levels. We observe that the mass drift decreases with mesh refinement, indicating that the injection-based scheme becomes more accurate as the mesh is refined. In contrast, the proposed conservative coarsening operator ensures exact mass conservation at each coarsening step, preventing mass drift over time.
\tikzexternaldisable
\tikzexternalenable
\subsection{Cahn--Hilliard}
\label{subsec:results_ch}

We begin with the Cahn--Hilliard equation, a stiff fourth-order nonlinear PDE that models phase separation in binary mixtures:
\begin{equation}
	\frac{\partial \phi}{\partial t} = \nabla \cdot \bigg(M (\phi) \nabla \bigg(\frac{\partial f}{\partial \phi} - \epsilon^2 \nabla^2 \phi\bigg)\bigg),
    \label{eq:CH_strong}
\end{equation}
where $\phi(\vec{x},t)$ is the concentration field, $M(\phi)$ is the mobility, $f(\phi)$ is the free energy of mixing, and $\epsilon$ controls the diffuse-interface thickness.

We focus on spinodal decomposition, in which an initially well-mixed state separates spontaneously into regions rich in each component. To quantify the thermodynamic state, we use the Ginzburg--Landau energy functional
\begin{equation}
E(\phi) = \int_{\Omega} \left( f(\phi) + \frac{1}{2}\epsilon^2 |\nabla \phi|^2 \right)\, d\Omega,
\label{eq:totalEnergy}
\end{equation}
where the first term is the bulk free energy and the second is the interfacial (gradient) contribution.

The CH dynamics exhibit two key properties that we use to assess numerical behavior under AMR:
\begin{enumerate}
\item \textit{Conservation of total mass:} $\displaystyle \int_{\Omega} \phi \, d\Omega$ is constant in time.
\item \textit{Monotonic decay of energy:} $\displaystyle E(\phi)$ decreases monotonically in time.
\end{enumerate}

\textbf{Free energy models:}
We consider (i) a polynomial double-well free energy
\begin{equation}
f(\phi) = \frac{1}{4} (1 - \phi^2)^2,
\label{eq:polynomial}
\end{equation}
for which $\phi\in[-1,1]$, and (ii) a logarithmic Flory--Huggins free energy (described in \secref{sec:FH_results}).

\textbf{Numerical discretization:}
To avoid a direct discretization of the fourth-order operator, we use the classical splitting strategy \citep{elliott1989second} by introducing the chemical potential $\mu$:
\begin{equation}
\begin{split}
\frac{\partial \phi}{\partial t} &= \nabla \cdot \big(M \nabla \mu\big), \\
\mu &= \frac{\partial f}{\partial \phi} - \epsilon^2 \nabla^2 \phi.
\end{split}
\label{eq:splitCH}
\end{equation}
We use continuous Galerkin finite elements with linear Lagrange basis functions for spatial discretization and backward Euler for time integration.

\textbf{Octree AMR strategy:}
Accurately resolving the diffuse interface is essential for capturing the correct phase separation dynamics. We therefore employ a two-level octree AMR strategy consisting of a fine \emph{interface} resolution $\mathcal{R}_i$ and a coarse \emph{bulk} resolution $\mathcal{R}_b$ ($\mathcal{R}_i>\mathcal{R}_b$).\footnote{A level-$l$ octree element corresponds to $\Delta x=\Delta y=\Delta z=L/2^l$, where $L$ is the domain length.} Elements satisfying $|\phi|\le \delta$ are marked as interface and maintained at the finest level $\mathcal{R}_i$, where $\delta$ is a user parameter controlling the interface band width.\footnote{For instance, $\delta=1$ refines the entire domain to $\mathcal{R}_i$, while $\delta=0$ keeps the domain at $\mathcal{R}_b$.} Adaptation is performed with balanced, single-level refinement/coarsening on the octree.

All simulations are performed with the massively parallel framework \texttt{Dendrite-kt}, which uses \texttt{Dendro-kt} for octree mesh generation and \texttt{PETSc} for linear algebra. Implementation details, scaling analyses, and validation for several multiphysics applications are reported in \citep{saurabh2021industrial,saurabh2022scalable,ishii2019solving,saurabh2021scalable}. The intergrid transfer operators described in \citet{saurabh2022scalable} form the baseline for the injection-based scheme and provide the foundation for the conservative coarsening operator proposed in \secref{sec:Algorithm}.

\subsubsection{2D spinodal decomposition (polynomial free energy) with linear basis functions}
\label{sec:2D_CH}

We consider $\Omega=[0,1]^2$ with an initial condition $\phi(\vec{x},0)=\phi_0+\eta(\vec{x})$, where $\phi_0=0$ and $\eta$ is a random uniform perturbation with zero mean and maximum amplitude $0.1$. We set $\epsilon^2=0.001$ and use constant mobility $M=1$. The timestep is $\Delta t=5\times 10^{-4}$. We choose bulk refinement $\mathcal{R}_b=4$ and interface refinement $\mathcal{R}_i=8$, and vary $\delta\in[0.90,0.98]$ to study sensitivity to the interface band width.

\figref{fig:evolution} shows the interface evolution: the initially mixed state separates into two phases; as time progresses, one phase forms islands that subsequently pinch-off and merge. \figref{fig:level} shows that the octree mesh resolves the interface at level $\mathcal{R}_i=8$ while retaining the bulk at level $\mathcal{R}_b=4$.

\figref{fig:CH_2D} compares conservation and energy behavior between coarsening schemes. We quantify mass drift by
\begin{equation}
	\Delta m(t) = \int_{\Omega} \phi(\vec{x},t)\, d\Omega \;-\; \int_{\Omega} \phi(\vec{x},0)\, d\Omega .
    \label{eq:mass_drift_def}
\end{equation}
As shown in \figref{fig:Mass_2D_CH}, the proposed scheme preserves mass to numerical precision for all $\delta$, while injection produces a nonzero drift that accumulates over time. The total energy $E(\phi)$ exhibits nearly identical monotonic decay for both schemes (\figref{fig:Energy_2D_CH}), indicating that the conservative coarsening does not degrade the physical dissipation behavior.

To isolate the effect of coarsening transfer on energetic consistency, we also compute the absolute energy mismatch across each coarsening event:
\begin{equation}
	\Delta E(t) = \left|E(\phi)_{\mathcal{M}_1} - E(\phi)_{\mathcal{M}_2}\right|,
    \label{eq:deltaE_def}
\end{equation}
where $\mathcal{M}_1$ and $\mathcal{M}_2$ denote the meshes immediately before and after coarsening, respectively. \figref{fig:Energy_diff_2D_CH} shows that the proposed scheme yields consistently smaller $\Delta E$ than injection by almost two orders of magnitude, reflecting improved fidelity of the transferred field under coarsening. 

We further compute the number of elements for both schemes over time (\figref{fig:Elements_2D_CH}) and observe almost a similar number of elements, indicating that differences in conservation and energy behavior are not due to differing mesh sizes. Higher value of $\delta$ leads to more conservative coarsening and thus have comparatively larger number of elements, irrespective of the coarsening scheme. Overall, both schemes yield similar mesh adaptation patterns, but the proposed conservative coarsening operator better preserves mass and energy consistency compared to the standard injection based approach.

\input{tikz/CH2D}

\begin{figure*}
	\centering
	\begin{subfigure}{0.3\textwidth}
		\begin{tikzpicture}[scale=0.7]
			\begin{axis}[
				xlabel={$t$},
				ylabel={$\Delta m$},
				legend pos=north west,
				xmin = 0,
				xmax = 3,
				]
				\addplot[mark=none, blue, very thick] table [col sep=comma, x=Time, y=massDiff_inj]{Data/2DCH/quadratic.csv};
				\addlegendentry{Injection};
				\addplot[mark=none, red, dashed, very thick] table [col sep=comma, x=Time, y=massDiff_consv]{Data/2DCH/quadratic.csv};
				\addlegendentry{Conservative};
			\end{axis}
		\end{tikzpicture}
		\caption{Mass difference plot}
		\label{fig:2D_CH_quad_Mass}
	\end{subfigure}
	\hspace{0.5cm}
	\begin{subfigure}{0.3\textwidth}
		\begin{tikzpicture}[scale=0.7]
			\begin{semilogxaxis}[
				xlabel={$t$},
				ylabel={$E(\phi)$},
				xmin = 5E-4,
				xmax = 3,
				yticklabel style={
					/pgf/number format/fixed,
				},
				scaled y ticks=false,
				]
				\addplot[mark=none, blue, very thick, each nth point=1] table [col sep=comma, x=Time, y=energy_injec] {Data/2DCH/quadratic.csv};
				\addlegendentry{Injection};
				\addplot[mark=none, red, dashed, very thick, each nth point=1] table [col sep=comma, x=Time, y=energy_consv] {Data/2DCH/quadratic.csv};
				\addlegendentry{Conservative};
			\end{semilogxaxis}
		\end{tikzpicture}
		\caption{Energy evolution}
		\label{fig:2D_CH_quad_Energy}
	\end{subfigure}
	\hspace{0.8cm}
	\begin{subfigure}{0.3\textwidth}
		\begin{tikzpicture}[scale=0.7]
			\begin{semilogyaxis}[		
				xmin = 0,
				xmax = 3,
				xlabel={$t$},
                legend pos = south west
				]
				\addplot[mark=none,very thick, blue]  table [col sep=comma, each nth point=4,x=Time, y=inj_diff] {Data/2DCH/quadratic.csv};
				\addlegendentry{Injection};
				\addplot[mark=none,very thick, red] table [ col sep=comma, each nth point=4,x=Time, y=consv_diff] {Data/2DCH/quadratic.csv};
				\addlegendentry{Conservative};
			\end{semilogyaxis}
		\end{tikzpicture}
		\caption{Energy difference}
		\label{fig:2D_CH_quad_EnergyDiff}
	\end{subfigure}
	\caption{Comparison of the energy and mass conservation between different interpolation schemes for quadratic basis function and polynomial free energy.}
	\label{fig:2D_CH_quad}
\end{figure*}
\subsubsection{2D spinodal decomposition (polynomial free energy) with quadratic basis function}
We further test the proposed scheme with quadratic basis functions for spatial discretization. We consider the same setup as in \secref{sec:2D_CH}, with $\Omega=[0,1]^2$, $\epsilon^2=0.001$, $M=1$, $\Delta t=5\times 10^{-4}$, $\mathcal{R}_b=3$, $\mathcal{R}_i=7$, and $\delta = 0.90$. \figref{fig:2D_CH_quad} shows results consistent with the linear basis function case. The proposed scheme preserves mass to numerical precision (\figref{fig:2D_CH_quad_Mass}), while injection exhibits significant mass drift. The total energy decay is similar across schemes (\figref{fig:2D_CH_quad_Energy}), but the coarsening-induced energy mismatch $\Delta E$ is substantially smaller for the proposed method (\figref{fig:2D_CH_quad_EnergyDiff}), indicating more accurate transfer during coarsening.

\subsubsection{3D spinodal decomposition (polynomial free energy)}
\label{sec:3D_CH}

\begin{figure*}
	\begin{subfigure}{0.32\textwidth}
		{\includegraphics[trim=300 20 500 20, clip,width=\linewidth]{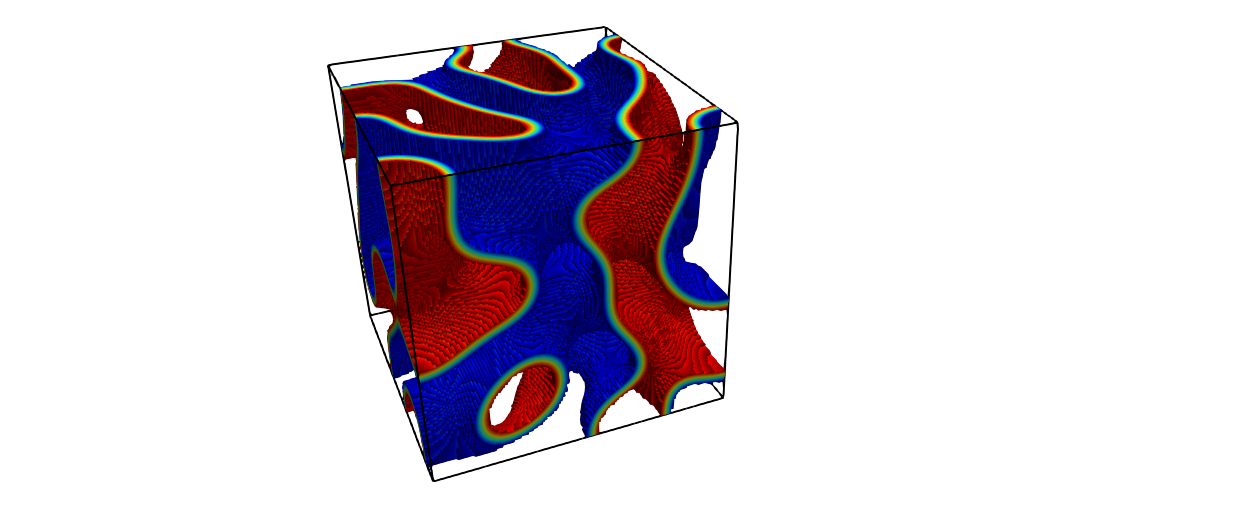}}
		\caption{$t = 0.1$}
		\label{fig:t0p1_3D}
	\end{subfigure}
	\begin{subfigure}{0.32\textwidth}
		{\includegraphics[trim=300 20 500 20, clip,width=\linewidth]{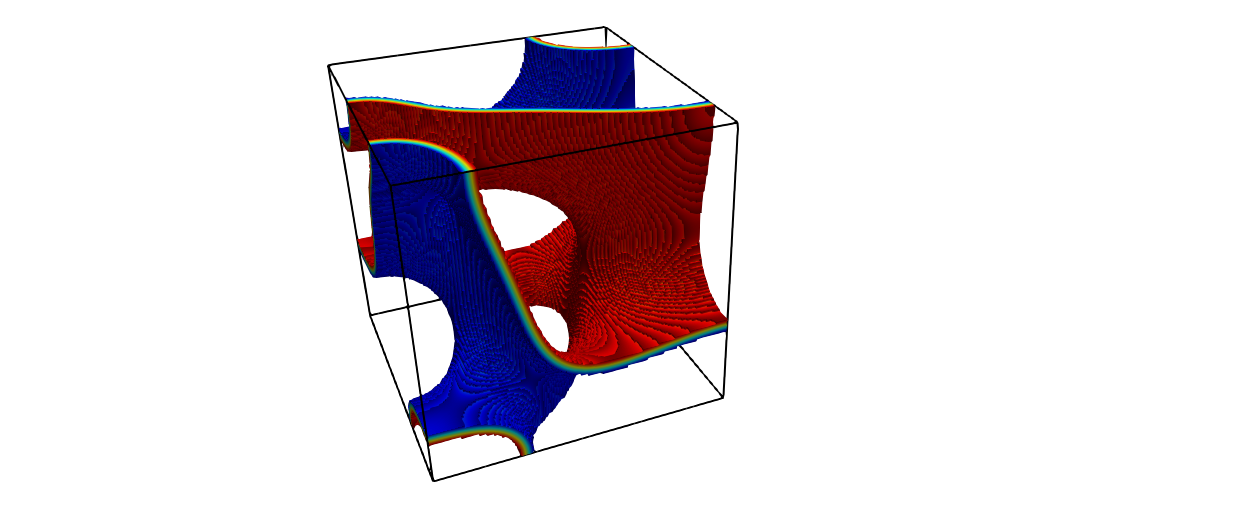}}
		\caption{$t = 1.0$}
		\label{fig:t1p0_3D}
	\end{subfigure}
	\begin{subfigure}{0.32\textwidth}
		{\includegraphics[trim=300 20 500 20, clip,width=\linewidth]{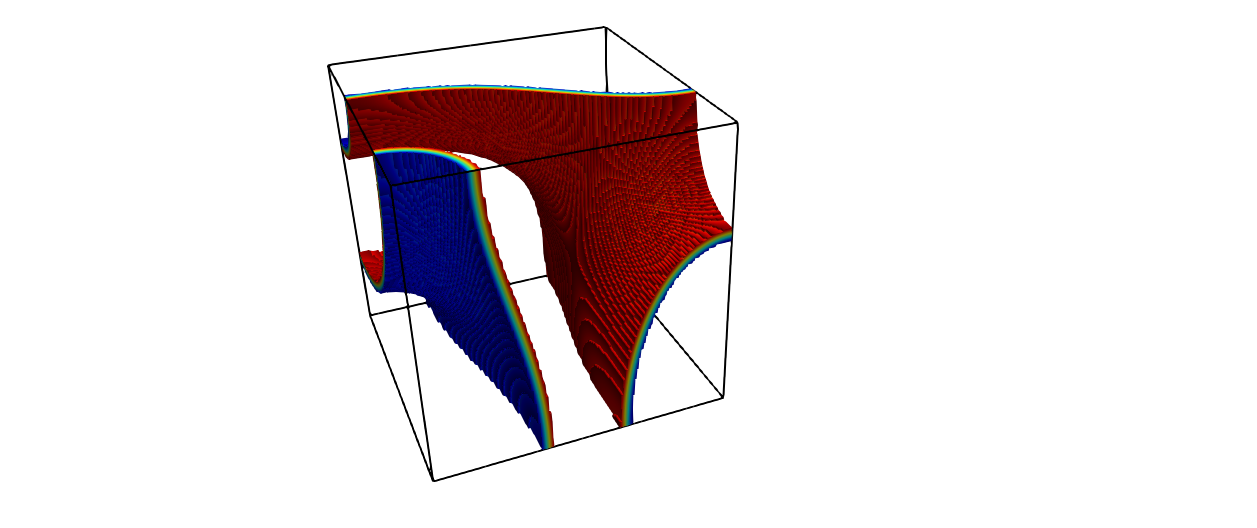}}
		\caption{$t = 3.0$}
		\label{fig:t3p0_3D}
	\end{subfigure}
	\caption{\textbf{3D case:} Evolution of the phase field variable $\phi$ over time for a uniform random mixture with perturbation around $\phi_0$ = 0.}
	\label{fig:evolution_3D}
\end{figure*}

We extend \secref{sec:2D_CH} to $\Omega=[0,1]^3$ using the same parameter values, except that the finest interface refinement is set to $\mathcal{R}_i=7$ (with $\mathcal{R}_b=4$) to control total DOF count. We consider $\delta=0.90$. \figref{fig:evolution_3D} shows interface evolution: small islands form at early times (\figref{fig:t0p1_3D}) and subsequently coarsen into larger domains (\figref{fig:t3p0_3D}).

\begin{figure*}
	\centering
	\begin{subfigure}{0.3\textwidth}
		\begin{tikzpicture}[scale=0.7]
			\begin{axis}[
				    xlabel={$t$},
				ylabel={$\Delta m$},
				legend pos=south east,
				      xmin = 0,
				xmax = 3,
				]
			 \addplot[mark=none, blue, very thick] table [col sep=comma, x=Time, y=massDiff_inj]{Data/3DCH/Polynomial_3D.csv};
				\addlegendentry{Injection};
				 \addplot[mark=none, red, dashed, very thick] table [col sep=comma, x=Time, y=massDiff_consv]{Data/3DCH/Polynomial_3D.csv};
				\addlegendentry{Conservative};
			\end{axis}
		\end{tikzpicture}
		\caption{Mass difference plot}
		\label{fig:3D_CH_mass}
	\end{subfigure}
\hspace{0.5cm}
	\begin{subfigure}{0.3\textwidth}
		\begin{tikzpicture}[scale=0.7]
			\begin{semilogxaxis}[
				    xlabel={$t$},
			    ylabel={$E(\phi)$},
			          xmin = 1E-3,
			    xmax = 3,
			     yticklabel style={
			    	/pgf/number format/fixed,
			    },
			    scaled y ticks=false,
				]
			    \addplot[mark=none, blue, very thick] table [col sep=comma, x=Time, y=energy_injec] {Data/3DCH/Polynomial_3D.csv};
			\addlegendentry{Injection};
			 \addplot[mark=none, red, dashed, very thick] table [col sep=comma, x=Time, y=energy_consv] {Data/3DCH/Polynomial_3D.csv};
			\addlegendentry{Conservative};
			\end{semilogxaxis}
		\end{tikzpicture}
		\caption{Energy evolution}
		\label{fig:3D_CH_Energy}
	\end{subfigure}
	\hspace{0.8cm}
		\begin{subfigure}{0.3\textwidth}
				\begin{tikzpicture}[scale=0.7]
						\begin{semilogyaxis}[		
						xmin = 1E-3,
						xmax = 3,
						xlabel={$t$},
						ylabel={$\Delta E(\phi)$},
						ymin=1E-13,
								]
								\addplot table [mark=none, col sep=comma, x=Time, y=inj_diff, each nth point=1] {Data/3DCH/Polynomial_3D.csv};
								\addlegendentry{Injection};
								\addplot table [mark=none, col sep=comma, x=Time, y=consv_diff, each nth point=1] {Data/3DCH/Polynomial_3D.csv};
								\addlegendentry{Conservative};
							\end{semilogyaxis}
					\end{tikzpicture}
					\caption{Energy difference}
					\label{fig:3D_CH_EnergyDiff}
			\end{subfigure}
	\caption{Comparison of the interpolation scheme for 3D Cahn-Hillard with polynomial free energy.}
	\label{fig:3D_CH}
\end{figure*}
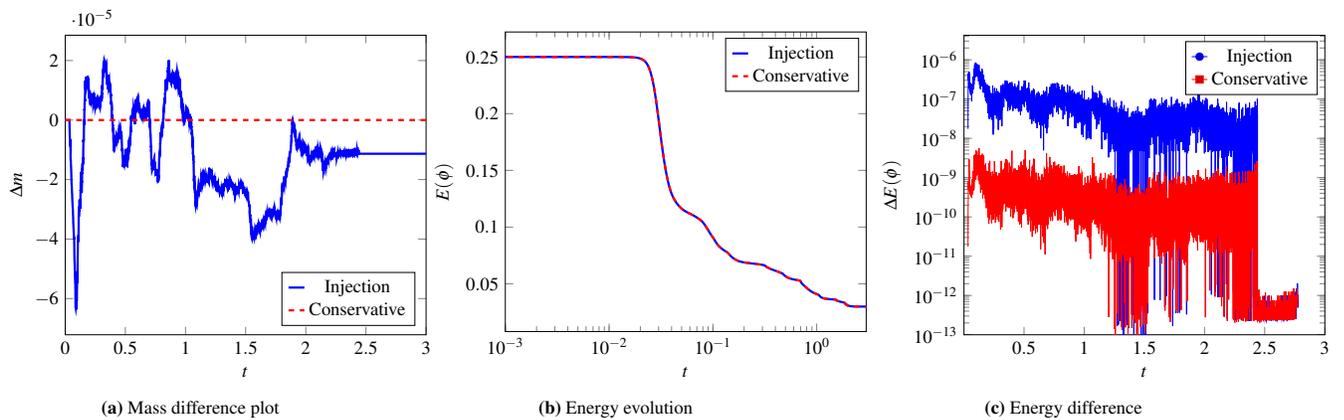

\figref{fig:3D_CH} shows trends consistent with the 2D case. The proposed scheme preserves mass to numerical precision (\figref{fig:3D_CH_mass}), whereas injection exhibits substantial mass loss that accumulates over time. The global energy decay remains similar across schemes, but the coarsening-induced energy mismatch $\Delta E$ is significantly smaller for the proposed method (\figref{fig:3D_CH_EnergyDiff}), indicating a more faithful transfer during coarsening.

\subsubsection{Flory--Huggins free energy}
\label{sec:FH_results}

We next test robustness to the free-energy model by using the logarithmic Flory--Huggins form
\begin{equation}
	f(\phi) = A \left(\phi \log(\phi) + (1 - \phi) \log(1 - \phi)\right)
    + \chi_{12}\,\phi (1 - \phi)
    + \beta \left(\frac{1}{\phi} + \frac{1}{1 - \phi}\right),
	\label{eq:FH}
\end{equation}
where $\chi_{12}$ is the Flory--Huggins interaction parameter and $\beta$ is a relaxation parameter. Here $\phi\in[0,1]$, with $\phi=0$ and $\phi=1$ representing the pure phases (in contrast to \eqnref{eq:polynomial}).

We use $A=1$, $\chi_{12}=3.0$, $\beta=0.01$, and $\epsilon^2=10^{-3}$. The initial condition is a random mixture with mean $0.5$ and fluctuation amplitude $0.01$. We take $\Delta t=5\times10^{-4}$ and define the interface band as $\delta\in(0.30,0.70)$. As before, $\mathcal{R}_b=4$, with $\mathcal{R}_i=8$ in 2D and $\mathcal{R}_i=7$ in 3D.

\figref{fig:2D_CH_FH} and \figref{fig:3D_CH_FH} show results for 2D and 3D, respectively. The same qualitative conclusions hold as in the polynomial case: the proposed scheme preserves mass to numerical precision while injection introduces noticeable mass fluctuations/drift (\figref{fig:2D_CH_FH_Mass}, \figref{fig:3D_CH_FH_Mass}). The total energy decay remains similar between schemes (\figref{fig:2D_CH_FH_Energy}, \figref{fig:3D_CH_FH_Energy}), but the proposed method produces substantially smaller coarsening-induced energy mismatch $\Delta E$ (\figref{fig:2D_CH_FH_EnergyDiff}, \figref{fig:3D_CH_FH_EnergyDiff}). These results demonstrate that the conservative coarsening operator is robust across different free-energy models.

\begin{figure*}
	\centering
	\hspace{-0.5cm}
	\begin{subfigure}{0.3\textwidth}
		\begin{tikzpicture}[scale=0.7]
			\begin{axis}[
				xlabel={$t$},
				ylabel={$\Delta m$},
				legend pos=north east,
				xmin = 0,
				xmax = 3,
				]
				\addplot[mark=none, blue, very thick] table [col sep=comma, x=Time, y=massDiff_inj]{Data/FH/FH_2D.csv};
				\addlegendentry{Injection};
				\addplot[mark=none, red, dashed, very thick] table [col sep=comma, x=Time, y=massDiff_consv]{Data/FH/FH_2D.csv};
				\addlegendentry{Conservative};
			\end{axis}
		\end{tikzpicture}
		\caption{Mass difference plot}
		\label{fig:2D_CH_FH_Mass}
	\end{subfigure}
	\hspace{0.7cm}
	\begin{subfigure}{0.3\textwidth}
		\begin{tikzpicture}[scale=0.7]
			\begin{semilogxaxis}[
				xlabel={$t$},
				ylabel={$E(\phi)$},
				xmin = 5E-4,
				xmax = 3,
				yticklabel style={
					/pgf/number format/fixed,
				},
				scaled y ticks=false,
				]
				\addplot[mark=none, blue, very thick, each nth point=1] table [col sep=comma, x=Time, y=energy_injec] {Data/FH/FH_2D.csv};
				\addlegendentry{Injection};
				\addplot[mark=none, red, dashed, very thick, each nth point=1] table [col sep=comma, x=Time, y=energy_consv] {Data/FH/FH_2D.csv};
				\addlegendentry{Conservative};
			\end{semilogxaxis}
		\end{tikzpicture}
		\caption{Energy evolution}
		\label{fig:2D_CH_FH_Energy}
	\end{subfigure}
	\hspace{0.8cm}
	\begin{subfigure}{0.3\textwidth}
		\begin{tikzpicture}[scale=0.7]
			\begin{semilogyaxis}[		
				xmin = 0,
				xmax = 3,
				xlabel={$t$},
				ylabel={$\Delta E(\phi)$},
				ymin=0,
				]
				\addplot[mark=none,very thick, blue]  table [col sep=comma, each nth point=4,x=Time, y=inj_diff] {Data/FH/FH_2D.csv};
				\addlegendentry{Injection};
				\addplot[mark=none,very thick, red] table [ col sep=comma, each nth point=4,x=Time, y=consv_diff] {Data/FH/FH_2D.csv};
				\addlegendentry{Conservative};
			\end{semilogyaxis}
		\end{tikzpicture}
		\caption{Energy difference}
		\label{fig:2D_CH_FH_EnergyDiff}
	\end{subfigure}
	\caption{Comparison of the interpolation scheme for 2D Cahn-Hillard with Flory Huggins free energy.}
	\label{fig:2D_CH_FH}
\end{figure*}

\begin{figure*}
	\centering
	\begin{subfigure}{0.3\textwidth}
		\begin{tikzpicture}[scale=0.7]
			\begin{axis}[
				xlabel={$t$},
				ylabel={$\Delta m$},
				legend pos=south east,
				xmin = 05E-4,
				xmax = 3,
				]
				\addplot[mark=none, blue, very thick, each nth point=4] table [col sep=comma, x=Time, y=massDiff_inj]{Data/FH/FH_3D.csv};
				\addlegendentry{Injection};
				\addplot[mark=none, red, dashed, very thick, each nth point=4] table [col sep=comma, x=Time, y=massDiff_consv]{Data/FH/FH_3D.csv};
				\addlegendentry{Conservative};
			\end{axis}
		\end{tikzpicture}
		\caption{Mass difference plot}
		\label{fig:3D_CH_FH_Mass}
	\end{subfigure}
	\hspace{0.5cm}
	\begin{subfigure}{0.3\textwidth}
		\begin{tikzpicture}[scale=0.7]
			\begin{semilogxaxis}[
				xlabel={$t$},
				ylabel={$E(\phi)$},
				xmin = 5E-4,
				xmax = 3,
				yticklabel style={
					/pgf/number format/fixed,
				},
				scaled y ticks=false,
				]
				\addplot[mark=none, blue, very thick] table [col sep=comma, x=Time, y=energy_injec] {Data/FH/FH_3D.csv};
				\addlegendentry{Injection};
				\addplot[mark=none, red, dashed, very thick] table [col sep=comma, x=Time, y=energy_consv] {Data/FH/FH_3D.csv};
				\addlegendentry{Conservative};
			\end{semilogxaxis}
		\end{tikzpicture}
		\caption{Energy evolution}
		\label{fig:3D_CH_FH_Energy}
	\end{subfigure}
	\hspace{0.8cm}
	\begin{subfigure}{0.3\textwidth}
		\begin{tikzpicture}[scale=0.7]
			\begin{semilogyaxis}[		
				xmin = 0,
				xmax = 3,
				xlabel={$t$},
				ylabel={$\Delta E(\phi)$},
				ymin=1E-11,
				]
				\addplot[mark=none,very thick, blue] table [ col sep=comma, x=Time, y=inj_diff] {Data/FH/FH_3D.csv};
				\addlegendentry{Injection};
				\addplot [mark=none,very thick,red] table [col sep=comma, x=Time, y=consv_diff] {Data/FH/FH_3D.csv};
				\addlegendentry{Conservative};
			\end{semilogyaxis}
		\end{tikzpicture}
		\caption{Energy difference}
		\label{fig:3D_CH_FH_EnergyDiff}
	\end{subfigure}
	\caption{Comparison of the interpolation scheme for 3D Cahn-Hillard with Flory Huggins free energy.}
	\label{fig:3D_CH_FH}
\end{figure*}

 \subsection{Cahn-Hillard Navier Stokes (CHNS)}
  \label{sec:bubble_case2}

 We further extended our conservation scheme to simulate two-phase flow using the CHNS set of partial differential equations (PDEs). The CHNS system offers a computationally tractable model that effectively captures interfacial dynamics in two-phase fluid flows. By employing Cahn-Hilliard equations to track the interface, we benefit from several advantages, such as mass conservation, thermodynamic consistency, and a free-energy-based description of surface tension, supported by well-established theory in non-equilibrium thermodynamics \cite{jacqmin1996energy}. Similar to Cahn-Hilliard, the CHNS system exhibits the same distinctive characteristics, namely mass conservation ($\int_{\Omega} \phi$ remains constant) and a monotonically decreasing total energy.
 
\textbf{Governing equations:}
 Let $\phi$ be the phase field variable that tracks the location of the phases and varies
 smoothly between +1 to -1 with a characteristic diffuse interface thickness $\epsilon$, and $v_i$ be the $i^{th}$ component of the mixture velocity of the two phases. $\eta_+$ and $\rho_+$ ($\eta_-$ and \added{$\rho_-$}) represent the viscosity and density of the +1 phase fluid (-1 phase fluid). The thermodynamically consistent CHNS equations~\citep{Khanwale2023projection, shen2015decoupled,anderson1998diffuse,feng2006fully} are written as follows:
 \begin{enumerate}[label=\roman*]
 	\item Momentum equation:
 	\footnotesize
 	\begin{equation*}
 		\begin{split}
 			&\pd{\left(\rho(\phi) v_i\right)}{t} + \pd{\left(\rho(\phi)v_iv_j\right)}{x_j} + \frac{1}{Pe}\pd{\left(J_jv_i\right)}{x_j}+\frac{Cn}{We} \pd{}{x_j}\left({\pd{\phi}{x_i}\pd{\phi}{x_j}}\right)
 			\\
 			&+\frac{1}{We}\pd{p}{x_i} - \frac{1}{Re}\pd{}{x_j}\left({\eta(\phi)\pd{v_i}{x_j}}\right) - \frac{\rho(\phi)\hat{{g_i}}}{Fr} = 0,\\
 			& \quad\quad \text{where,} \quad J_i = \frac{\left(\rho_- - \rho_+\right)}{2\;\rho_+ Cn} \, m(\phi)\pd{\mu}{x_i},
 		\end{split}
 	\end{equation*}
 	\normalsize
 	\item Solenoidality and Continuity:
 	\footnotesize
 	\begin{equation*}
 		\pd{v_i}{x_i} = 0, \quad \pd{\rho(\phi)}{t} + \pd{\left(\rho(\phi)v_i\right)}{x_i}+
 		\frac{1}{Pe} \pd{J_i}{x_i} = 0,
 	\end{equation*}
 	\normalsize
 	\item Cahn--Hillard Equations
 	\footnotesize
 	\begin{equation*}
 		\begin{split}
 			&\pd{\phi}{t} + \pd{\left(v_i \phi\right)}{x_i} - \frac{1}{PeCn} \pd{}{x_i}\left(m(\phi){\pd{\mu}{x_i}}\right) = 0 \\
 			&\mu = \psi'(\phi) - Cn^2 \pd{}{x_i}\left({\pd{\phi}{x_i}}\right)
 		\end{split}
 	\end{equation*}
 	\normalsize
 \end{enumerate}
 
 where $\rho(\phi)$ is the non--dimensional mixture density given by $\left({\rho_+ - \;\rho_-}/{2\rho_+}\right) \phi + \left({\rho_+ + \;\rho_-}/{2\rho_+}\right)$, the non--dimensional mixture viscosity $\eta\left(\phi\right)$ is given by $\left({\eta_+ - \;\eta_-}/{2\eta_+}\right) \phi + \left({\eta_+ + \;\eta_-}/{2\eta_+}\right)$.

 Non-dimensional parameters are as follows: Peclet, $Pe = \frac{u_{r} L_{r}^2}{m_{r}\sigma}$; Reynolds, $Re = \frac{u_{r} L_{r}}{\nu_{r}}$; Weber, $We = \frac{\rho_{r}u_{r}^2 L_{r}}{\sigma}$; Cahn, $Cn = \frac{\varepsilon}{L_{r}}$; and Froude, $Fr = \frac{u_{r}^2}{gL_{r}}$, with $u_{r}$ and $L_r$ denoting the reference velocity and length, respectively. The mobility $m(\phi)$ is set to a constant value of 1.

 Several studies have focused on developing robust numerical methods that uphold the aforementioned properties of the CHNS system. However, previous results have demonstrated a discrepancy in mass conservation when these methods are combined with adaptive mesh refinement (AMR). 
 We employ the projection-based discretization scheme proposed by \citet{Khanwale2023projection} to simulate the system. Our simulation utilizes second-order semi-implicit time stepping in conjunction with linear basis functions to achieve second-order spatial accuracy.
 
 We consider a 2D bubble rise as a benchmark problem~\citep{hysing2009quantitative,aland2012benchmark,khanwale2022fully,Khanwale2023projection} and compare the result with our proposed mass-conserving interpolation scheme. \tabref{tab:bubble-rise} shows the parameter used for the simulation. Similar to the Cahn-Hillard case, the AMR procedure is defined by two levels of refinement: bulk ($\mathcal{R}_b$) and interface ($\mathcal{R}_i$), which is set to 6 and 11 respectively. We chose the value of $\delta$ to be 0.95 for the refinement threshold. As the interface moves, the element close to the interface is refined and coarsened away from it. 
 
\begin{table*}[htb]
	\centering
	\begin{tabular}{|c|c|c|c|c|c|c|c|c|c|c|c|c|c|c|c|}
		\hline
		$\rho_+$ & $\rho_-$ & $\mu_+$ & $\mu_-$ & $\rho+/\rho-$ & $\nu_+ / \nu_-$ & $g$ & $\sigma$ & $Ar$ & $We$ & $Fr$ & $Cn$ & $Pe$ & $\Delta t$ & $\mathcal{R}_b$ & $\mathcal{R}_i$ \\ \hline
		1000 & 1.0 & 10 & 0.1 & 1000 & 100 & 0.98 & 1.96 & 35 & 125 & 1.0 & 0.01 & 3333.33 & 0.0025 & 6 & 11 \\ \hline
	\end{tabular}
	\caption{Table showing the physical and simulation parameter for the bubble rise benchmark problem considered from \citet{Khanwale2023projection}.	}
	\label{tab:bubble-rise}
\end{table*}
%

 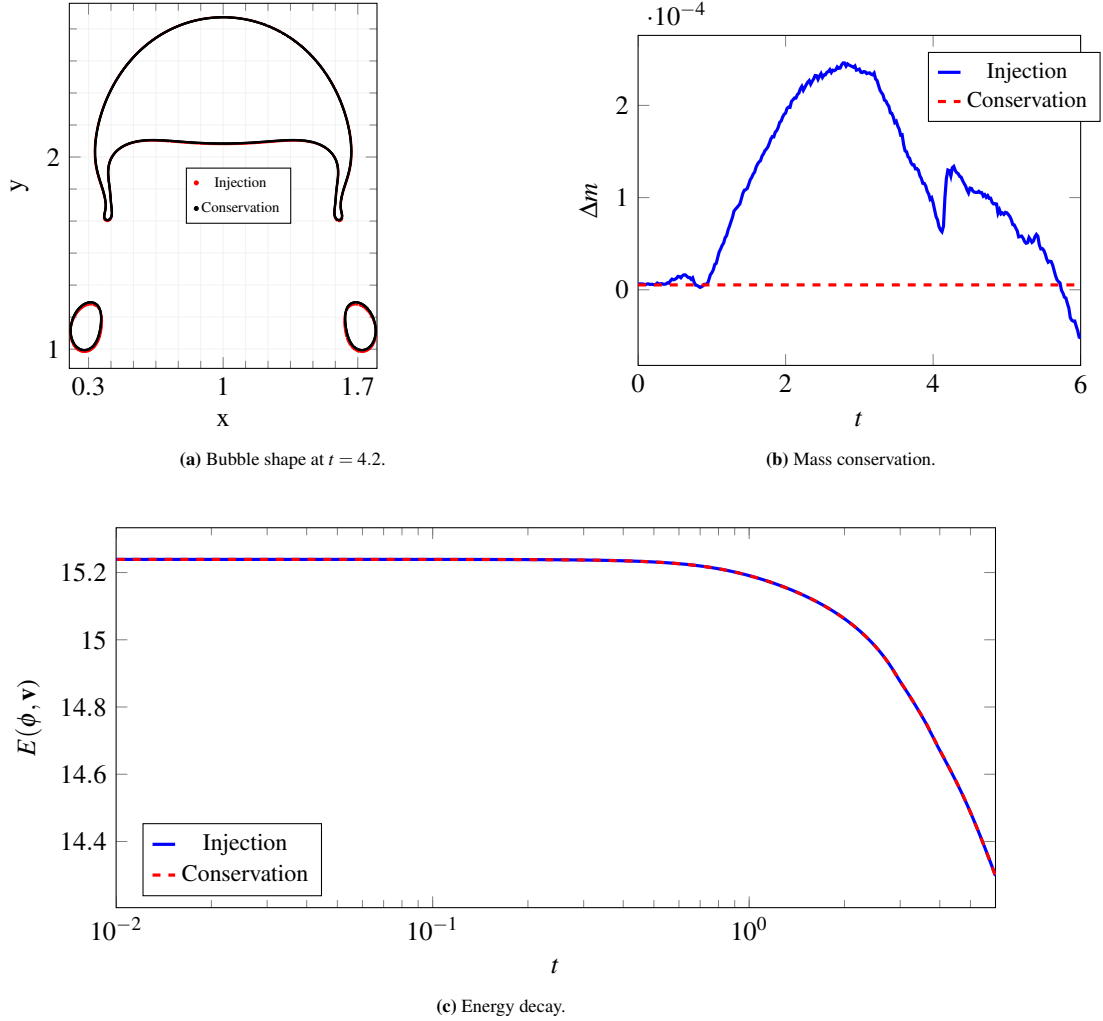
\begin{figure}[h]
\begin{subfigure}{0.45\linewidth}
	
    \begin{tikzpicture}
\begin{axis}[
          width=1.0\linewidth, 
          scaled y ticks=true,
          xlabel={$\mathrm{x}$},
          ylabel={$\mathrm{y}$},
legend style={nodes={scale=0.65, transform shape}}, ymin=0.9, ymax=2.8, ytick distance=1.0,  xtick={0.3, 1.0, 1.7},
legend style={at={(0.55,0.55)},anchor= north,
	nodes={scale=0.95, transform shape}, row sep=2.5pt},
	unit vector ratio*=1 1 1,
	xmin=0.2, xmax=1.8, 
	legend image post style={scale=3.0},
	grid=both,
	grid style={line width=.1pt, draw=gray!10},
	minor tick num=5
	]
]

    \addplot [only marks,mark size = 0.25pt,color=red, each nth point=8] table [x expr=\thisrow{"Points:0"},y expr=\thisrow{"Points:1"}, col sep=comma] {Data/BubbleRise2DCase2/phi0_Injection.csv};
    
    \addplot [only marks,mark size = 0.25pt,color=black, each nth point=4] table [x expr=\thisrow{"Points:0"},y expr=\thisrow{"Points:1"}, col sep=comma] {Data/BubbleRise2DCase2/phi0MConv.csv};

\legend{\small  Injection ,\small Conservation}

  \end{axis}
  \end{tikzpicture}
\caption{Bubble shape at $t = 4.2$.}
\label{fig:bubble2_shape}
  \end{subfigure}
  \begin{subfigure}{0.45\linewidth}
	\begin{tikzpicture}
		\begin{axis}[
			width=1.0\linewidth, 
			height = 0.8\linewidth,
			xmin = 0,
			xmax = 6,
			xlabel={$t$},
			ylabel = {$\Delta m$},
				legend style={at={(0.85,0.95)},anchor= north,nodes={scale=0.95, transform shape}}
			]
			\addplot[mark=none,color = blue, very thick,each nth point=10] 
			table[x expr={\thisrow{time}))},y expr=(\thisrow{TotalPhi}-6.429198800557858640e+00),col sep=space]{Data/BubbleRise2DCase2/InjectionEnergy_data.txt};
			
			\addplot[mark=none,color = red, very thick,each nth point=10,dashed] 
			table[x expr={\thisrow{time}))},y expr=(\thisrow{TotalPhi} - 6.429198800557858640e+00),col sep=space]{Data/BubbleRise2DCase2/MConvEnergy_data.txt};
			\legend{\small  Injection ,\small Conservation}
		\end{axis}
	\end{tikzpicture}
\caption{Mass conservation.}
\label{fig:bubble2_mass}
\end{subfigure}

\vspace{2em}
  \begin{subfigure}{0.8\linewidth}
  	\centering
    \begin{tikzpicture}
	\begin{semilogxaxis}[
		width=1.0\linewidth, 
		height = 0.5\linewidth,
		xmin = 1E-2,
		xmax = 6,
		xlabel={$t$},
		ylabel={$E(\phi,\Vector{v}$)},
		legend pos=south west
		]
	 \addplot[mark=none,color = blue, very thick,each nth point=2] 
	table[x expr={\thisrow{time}))},y expr=(\thisrow{TotalEnergy}),col sep=space]{Data/BubbleRise2DCase2/InjectionEnergy_data.txt};
	
	 \addplot[mark=none,color = red, very thick,each nth point=2,dashed] 
	table[x expr={\thisrow{time}))},y expr=(\thisrow{TotalEnergy}),col sep=space]{Data/BubbleRise2DCase2/MConvEnergy_data.txt};
	        \legend{\small  Injection ,\small Conservation}
	\end{semilogxaxis}
\end{tikzpicture}
\caption{Energy decay.}
\label{fig:bubble2_energy}
\end{subfigure}
\caption{2D Single Rising Drop Test Case 2: Figure showing the comparison for proposed mass conserving interpolation with injection-based interpolation for the test case in \secref{sec:bubble_case2}}
\label{fig:bubble_case2}
\end{figure}
 \figref{fig:bubble_case2} shows the result for the simulation. We note that only the variable $\phi$ is interpolated using mass conservation interpolation, whereas the other variable $(\vec{v}, p, \mu)$ is interpolated using the injection scheme during the coarsening stage. We do not observe any major difference in the bubble's shape from our proposed mass conservation scheme (\figref{fig:bubble2_shape}).
 However, the proposed scheme effectively preserves the total mass of the system. In contrast, the injection-based coarsening procedure leads to a mass discrepancy of approximately $\mathcal{O}(10^{-4})$ (\figref{fig:bubble2_mass}), while the proposed scheme demonstrates mass conservation, as indicated by the nearly flat line. Additionally, the total energy of the system, as defined by \eqnref{eqn:energy_functional}, shows a monotonic decay over time, with no significant difference between the two  interpolation schemes (\figref{fig:bubble2_energy}).
 
 \begin{equation}
 	\begin{split}
 		E_{\text{tot}}(\vec{v},\phi, t) =  & \int_{\Omega}\frac{1}{2}\rho(\phi) \norm{\vec{v}}^2 \mathrm{d}\vec{x} \\
 		& \frac{1}{CnWe}\int_{\Omega} \left(\psi(\phi) + \frac{Cn^2}{2} \norm{\nabla\phi}^2 + \frac{1}{Fr} \rho(\phi) y \right) \mathrm{d}\vec{x}.
 	\end{split}
 	\label{eqn:energy_functional}
 \end{equation}

\subsection{Computational cost and parallel scalability}
We implemented the overall procedure in the code framework \textsc{Dendrite-kt}. We have presented the scalability results till $\mathcal{O}(114K)$  of the overall framework in our previous works~\citep{saurabh2022scalable,saurabh2021industrial,khanwale2022breakupdynamicsprimaryjet, Khanwale2023projection,khanwale2020simulating,saurabh2021scalable}. The two additional steps that are required to perform the proposed scheme are:

\begin{enumerate}
	\item Perform the local $L^2$ projection from the fine mesh to the coarse mesh as described in \Algref{alg:l2_restrict_dd} to transfer the values from the fine elements to the coarse elements at the Gauss points.
	\item Invert the constant coefficient global mass matrix.
\end{enumerate}

The coarsening stage begins by constructing a surrogate tree from the original octree to ensure that all leaf elements of a parent marked for coarsening reside on the same processor. This redistribution step is required regardless of the intergrid transfer scheme and is described in detail in \citet{saurabh2022scalable}. This is done to ensure that the subsequent coarsening operations can be performed locally without inter-process communication. After establishing the surrogate tree, we execute the two steps outlined above.

The first part is completely local and independent for each coarsened element, involving only the children elements and their corresponding entries in the fine elements. The step do not introduce any additional inter-process communication. Additional details can be found in Section II.C of \citet{saurabh2022scalable}. This leverages the already existing algorithm on octrees and have been demonstrated to be scalable in our previous works. The second part relies on the optimal inversion of the mass matrix. This can be treated as a standard linear system inversion. Overall, the overhead of the proposed method is comparatively very small, which is mostly dominated by the cost of interface tracking, remeshing, and pressure projection solve for CHNS. 

In our experiments, we observe that the additional cost incurred by the proposed method is approximately 10\% of the total computational cost per time step for Cahn-Hillard solve, and about 7\% for the 2D CHNS simulation. This overhead is expected to decrease further for larger problem sizes, as the cost of the mass matrix inversion becomes less significant relative to the overall simulation which will be mostly dominated by the cost of pressure projection and interface tracking.

\section{Conclusion}

We presented a simple and scalable \emph{field-conserving coarsening} operator for adaptive mesh refinement (AMR) on massively parallel, balanced \emph{octree} meshes. The central observation is that, for continuous Galerkin (CG) discretizations, parent-to-child refinement transfer is naturally conservative, whereas conventional injection-based child-to-parent coarsening transfer is generally not. Under repeated AMR cycles, this nonconservative coarsening can introduce systematic drift (or oscillations) in globally conserved quantities, which becomes particularly problematic in long-horizon simulations.

To eliminate this failure mode, we proposed a conservative coarsening procedure that enforces discrete global conservation during coarsening by (i) computing conservative coarse-element values at quadrature points using a local $L^2$ projection at quadrature points and (ii) recovering coarse nodal degrees of freedom via a global $L^2$ projection, i.e., a mass-matrix solve on the coarse mesh. While this introduces an additional linear solve compared to injection, the operator is highly scalable on distributed-memory architectures and can be implemented efficiently within existing octree AMR frameworks.

We validated the method on mass-conserving phase-field models, including the Cahn--Hilliard equation and the coupled Cahn--Hilliard--Navier--Stokes (CHNS) system. Across 2D and 3D tests and multiple free-energy models, the proposed scheme preserved mass to numerical precision under AMR, while injection exhibited noticeable mass drift. Importantly, the conservative coarsening achieved this without degrading the expected energetic behavior, and it reduced coarsening-induced energy mismatch by almost two orders of magnitude relative to injection.

Although we emphasized mass conservation in this work, the same framework can be extended to enforce discrete conservation of other integral quantities (e.g., species, charge, or application-specific invariants) by modifying the conserved moment(s) imposed during coarsening \citep{baiges2021adaptive}. Future work will explore (i) enforcing multiple conserved quantities simultaneously, (ii) extensions to higher-order CG discretizations with Lagrange basis functions of arbitrary degree, and (iii) assessing the impact of conservation enforcement on accuracy, stability, and solver performance in broader multiphysics applications. Overall, the proposed conservative coarsening operator offers a practical approach to enhancing the robustness and reliability of octree-AMR simulations, particularly when long-term conservation is crucial.

\section{Acknowledgments}
BG and HS acknowledge funding from NSF 2513871. We acknowledge funding from NSF 2323716 and NSF 2053760.

\bibliographystyle{./bibliography/elsarticle-num-names}
\bibliography{./bibliography/main}
\appendix
\section{Restriction matrices for 1D elements}
\label{sec:Restriction_matrices}
In this appendix, we provide the restriction matrices for Q1 and Q2 elements in 1D used in \secref{sec:results} for constructing the coarse element quadrature values from the fine element quadrature values during the coarsening operation. We have used standard Gauss-Legendre quadrature for both Q1 and Q2 elemnts leading to $N_{\mathrm{ip},c} = N_{\mathrm{ip},f} = 2$ and $N_{\mathrm{ip},c} = N_{\mathrm{ip},f} = 3$ quadrature points per coarse and fine element, respectively. The resulting restriction matrices of the size $N_{\mathrm{ip},c}  \times 2 N_{\mathrm{ip},f}$  are given below.
\subsection{Restriction matrix for Q1 element in 1D}
The restriction matrix for Q1 element for standard Gauss Legendre quadrature with $N_{\mathrm{ip},c} = N_{\mathrm{ip},f} = 2$ in 1D is given by:

\begin{equation}
    \small
  R_{1D} = \begin{pmatrix}
0.5915063509461096 & 0.3415063509461096 & 0.1584936490538904 & -0.09150635094610965 \\
-0.09150635094610965 & 0.1584936490538904 & 0.3415063509461096 & 0.5915063509461096
\end{pmatrix}
\end{equation}

\subsection{Restriction matrix for Q2 element in 1D}
The restriction matrix for Q2 element for standard Gauss Legendre quadrature with $N_{\mathrm{ip},c} = N_{\mathrm{ip},f} = 3$ in 1D is given by:
\begin{equation}
    \small
  R_{1D} = \begin{pmatrix}
0.614415278851 & 0.424865556414 & 0.041666666667 & -0.031081945517 & -0.091532223080 & 0.041666666667 \\
-0.097551215948 & 0.291666666667 & 0.305884549282 & 0.305884549282 & 0.291666666667 & -0.097551215948 \\
0.041666666667 & -0.091532223080 & -0.031081945517 & 0.041666666667 & 0.424865556414 & 0.614415278851
\end{pmatrix}
\end{equation}

The above restriction matrices can be used in \eqnref{eq:restriction_matrix_form} to compute the coarse element quadrature values from the fine element quadrature values during the coarsening operation in 1D. The extension to higher dimensions can be achieved via tensor products of the 1D restriction matrices along each dimension (\Algref{alg:l2_restrict_dd}).
\end{document}